\documentclass[draft,reqno]{amsproc}
\usepackage{amsthm}
\usepackage{amsfonts}
\usepackage{amsmath}
\usepackage{amssymb}
\usepackage{mathrsfs}
\usepackage{euscript}   
\usepackage{color}
\usepackage{tikz}

\makeatletter
\@namedef{subjclassname@2010}{%
\textup{2010} Mathematics Subject Classification}
\makeatother

\numberwithin{equation}{section}

\newtheorem{thm}{Theorem}[section]
\newtheorem{pro}[thm]{Proposition}
\newtheorem{cor}[thm]{Corollary}
\newtheorem{lem}[thm]{Lemma}

\newtheorem*{opq*}{\bf Problem}
\newtheorem{dfn}[thm]{Definition}

\theoremstyle{remark}
\newtheorem{rem}[thm]{Remark}

\theoremstyle{definition}
\newtheorem{exa}[thm]{Example}

\newcommand*{\D}{\mathrm{d\hspace{.1ex}}}
\newcommand*{\alphab}{\boldsymbol{\alpha}}

\newcommand*{\borel}[1]{{\mathfrak B}(#1)}
\newcommand*{\cbb}{\mathbb{C}}
\newcommand*{\dbb}{\mathbb{D}}

\newcommand*{\lcal}{\EuScript{L}}
\newcommand*{\ccal}{\EuScript{C}}
\newcommand*{\gammab}{\boldsymbol{\gamma}}

\newcommand*{\Ge}{\geqslant}
\newcommand*{\hh}{\EuScript{H}}
\newcommand*{\mm}{\EuScript{M}}

\newcommand*{\I}{{\mathrm i\hspace{.1ex}}}

\newcommand*{\is}[2]{\langle#1,#2\rangle}
\newcommand*{\jd}[1]{\mathscr N(#1)}
\newcommand*{\kk}{\EuScript{K}}
\newcommand*{\Le}{\leqslant}

\newcommand*{\mcal}{\EuScript{M}}
\newcommand*{\nbold}{\boldsymbol{N}}
\newcommand*{\mbold}{\boldsymbol{M}}

\newcommand*{\nbb}{\mathbb N}
\newcommand*{\ob}[1]{{\mathscr R}(#1)}

\newcommand*{\ogr}[1]{\boldsymbol B(#1)}
\newcommand*{\rbb}{\mathbb R}
\newcommand{\restr}[2]{\left.#1\right|_{#2}}
\newcommand*{\sbd}{\bar{\mathbb{B}}_2}
\newcommand*{\sbb}{\bar{\mathbb{B}}^{+}_2}
\newcommand*{\sfrak}{\mathfrak{S}}
\newcommand*{\smat}[1]{\left[\begin{smallmatrix}#1\end{smallmatrix}\right]}

\newcommand*{\supp}{\mathrm{supp}\,}
\newcommand*{\tbb}{\mathbb{T}}
\newcommand*{\tbold}{\boldsymbol{T}}
\newcommand*{\sbold}{\boldsymbol{S}}
\newcommand*{\qbold}{\boldsymbol{Q}}
\newcommand*{\zbb}{\mathbb Z}
   
\begin{document}
   \title[
Subnormality of $B$-operators of class $\mathscr{S}$]
{Subnormality of $B$-operators of class $\mathscr{S}$
   }
   \author[S. Chavan]{Sameer Chavan}
   \address{Department of Mathematics and Statistics\\
Indian Institute of Technology Kanpur, India}
   \email{chavan@iitk.ac.in}
   \author[Z.\ J.\ Jab{\l}o\'nski]{Zenon Jan
Jab{\l}o\'nski}
   \address{Instytut Matematyki,
Uniwersytet Jagiello\'nski, ul.\ \L ojasiewicza 6,
PL-30348 Kra\-k\'ow, Poland}
\email{Zenon.Jablonski@im.uj.edu.pl}
   \author[I.\ B.\ Jung]{Il Bong Jung}
   \address{Department of Mathematics,
Kyungpook National University, Daegu 702-701, Korea}
   \email{ibjung@knu.ac.kr}
   \author[J.\ Stochel]{Jan Stochel}
\address{Instytut Matematyki, Uniwersytet
Jagiello\'nski, ul.\ \L ojasiewicza 6, PL-30348
Kra\-k\'ow, Poland} \email{Jan.Stochel@im.uj.edu.pl}
   \thanks{The research of the first author was supported by the
Advanced Research Grant ANRF/ARG/2025/001228MS. The
research of the remaining authors was supported by the
National Science Centre, Poland, under project OPUS
grant no. DEC-2024/55/B/ST1/00280.}
   \subjclass[2020]{Primary 47B20, 47A13; Secondary
47B28, 47A10, 44A60}
   \keywords{$B$-operator, (joint) normal extension,
(joint) quasinormal extension, Taylor spectrum,
spherical contractivity}
   \maketitle
   \begin{abstract}
We characterize subnormality of $B$-operators
$T=\bigl[\begin{smallmatrix}V&E\\0&S\end{smallmatrix}\bigr]$
on $\hh_1\oplus\hh_2$, where $V$ is an isometry,
$V^*E=0$, $S$ is subnormal, and $S|E|=|E|S$. We
determine when $T$ admits compactly supported signed
elementary measures and express them explicitly in
terms of the elementary measures of $S$. Requiring the
signed elementary measures of $T$ to be positive
yields a complete measure-theoretic characterization
of the subnormality of $T$. If $E$ is injective, $T$
is subnormal if and only if it is a contraction. In
the general case, the subnormality of $T$ is
equivalent to the contractivity of its restriction to
$\hh_1\oplus\overline{\ob{|E|}}$. We also obtain
equivalent criteria involving the Taylor spectra of
$(S,|E|)$, its minimal normal and minimal quasinormal
extensions, and the corresponding pairs of moduli.
These results extend the characterization previously
known for a quasinormal lower-right entry.
   \end{abstract}
   \section{Introduction}
Block operator matrices of size $2\times 2$ form an
important class of operators in functional analysis
and operator theory. Their structural properties make
them a powerful tool in the study of a wide range of
problems, including the Halmos similarity problem for
polynomially bounded operators \cite{Fog64,Pi97}, the
hyperinvariant subspace problem
\cite{Do-Pe72,Ki12,JKP18,JKP19}, the description of
the Taylor spectrum of pairs of Hilbert space
operators \cite{Vas78,Vas82}, model-theoretic
approaches to the time-shift operator in stochastic
process theory \cite{Ag-St}, criteria for the
invertibility of upper triangular block matrices
\cite{WLee00}, investigations of $2$-hyponormal
operators \cite{Cur-Lee02}, and the classification of
Cowen-Douglas operators via complete unitary
invariants derived from upper triangular $2\times 2$
representations \cite{Mis17}.

Brownian-type operators, or $B$-operators, arise from
the theory of Brownian unitaries developed by Agler
and Stankus \cite{Ag-St}. The study of subnormality in
this class was initiated in \cite{C-J-J-S21}, and a
related extension theory was developed in
\cite{C-J-J-S24}. In this paper we characterize
subnormality for $B$-operators of class $\mathscr{S}$,
in which the lower-right entry is subnormal. This
extends \cite[Theorem~1.2]{C-J-J-S21}, which treats
the case of a quasinormal entry.

For such an operator
$T=\bigl[\begin{smallmatrix}V&E\\0&S\end{smallmatrix}\bigr]$
acting on $\hh_1\oplus\hh_2$, the defining relations
yield (see Lemma~5.1)
   \begin{align*}
\|T^n(0,f)\|^2 =\|S^n
f\|^2+\sum_{j=0}^{n-1}\|S^j|E|f\|^2, \qquad
n\Ge1,\quad f\in\hh_2.
   \end{align*}
When $S$ is quasinormal, the identity
$S^{*n}S^n=|S|^{2n}$ allows the squared norms
$\|T^n(0,f)\|^2$ to be computed using the joint
functional calculus of $(|S|,|E|)$. For a general
subnormal $S$, this reduction is unavailable, and we
work instead with the elementary measures $\mu_f$
representing the sequences $\{\|S^n
f\|^2\}_{n=0}^{\infty}$.

Theorem~\ref{wkyv} determines exactly when the
sequences $\{\|T^n(0,f)\|^2\}_{n=0}^{\infty}$ admit
finite signed representing measures $\rho_f$ with
compact support and gives their explicit form. The
criterion is governed by the behavior at $1$ of the
elementary measure $\mu_{|E|f}$: an atom at $1$ is
excluded, and $|1-x|^{-1}$ must be integrable with
respect to $\mu_{|E|f}$. The positivity of $\rho_f$
for all $f\in\hh_2$ then yields a characterization of
subnormality by a measure domination condition on
$[0,1)$ and a vanishing condition on $[1,\infty)$
(Theorem~\ref{mnteo}).

When $E$ is injective, these measure conditions are
equivalent to the single operator inequality
$S^*S+E^*E\Le I$. Thus $T$ is subnormal if and only if
it is a contraction, or, equivalently, if and only if
$(S,|E|)$ is a spherical contraction
(Theorem~\ref{wkwsc}). An intermediate
characterization is given by a factorization of $|E|$
through the defect operator of a minimal normal
extension of $S$, with a contractive factor commuting
with the modulus of that extension
(Theorem~\ref{cvxw}). For general $E\neq0$, the
subnormality of $T$ is equivalent to the contractivity
of its restriction to $\hh_1\oplus\overline{\ob{|E|}}$
(Proposition~\ref{13posr}).

The resulting spectral characterization takes the form
   \begin{align*}
\sigma(S,|E|)\subseteq \sbd\cup(\cbb\times\{0\}),
   \end{align*}
where $\sigma$ denotes the Taylor spectrum and $\sbd$
the closed unit ball of $\cbb^2$. Equivalent criteria
hold for $(|S|,|E|)$, for minimal normal and minimal
quasinormal extensions of $(S,|E|)$, and for the
corresponding pairs of moduli (Theorem~\ref{aswq}).
When $E$ is injective, the spectral condition reduces
to containment in $\sbd$. The appendix establishes the
existence and uniqueness, up to unitary equivalence
fixing the original space, of minimal joint
quasinormal extensions of subnormal pairs.

We conclude this introduction with the notation and
basic definitions used throughout the paper. All
Hilbert spaces are complex. For Hilbert spaces $\hh$
and $\kk$, let $\ogr{\hh,\kk}$ denote the Banach space
of bounded linear operators from $\hh$ to $\kk$, and
write $\ogr{\hh}:=\ogr{\hh,\hh}$. The kernel, range,
adjoint, and modulus of an operator $T$ are denoted by
$\jd{T}$, $\ob{T}$, $T^*$, and $|T|$, respectively. We
use $I=I_\hh$ for the identity operator, $\sigma(T)$
and $\sigma_{\mathrm p}(T)$ for the spectrum and point
spectrum of $T\in\ogr{\hh}$, and $\ccal'$ for the
commutant of $\ccal\subseteq\ogr{\hh}$.

An operator $N\in\ogr{\kk}$ is \emph{normal} if
$N^*N=NN^*$. An operator $S\in\ogr{\hh}$ is
\emph{subnormal} if it is the restriction of a normal
operator $N\in\ogr{\kk}$ to an invariant subspace
$\hh\subseteq\kk$. Such a \emph{normal extension} is
\emph{minimal} if $\kk=\bigvee_{n\Ge0}N^{*n}\hh$;
every subnormal operator admits one
\cite[Proposition~II.2.4]{Co91}. For a minimal normal
extension $N$ of $S$, the Bram lifting theorem
\cite{Br55} gives a unital $*$-isomorphism

   \begin{align*}
\{S,S^*\}'\ni A\longmapsto\widehat A
\in\{N\}'\cap\{P\}',
   \end{align*}
where $P$ is the orthogonal projection of $\kk$ onto
$\hh$ and $\widehat A|_{\hh}=A$ (see also
\cite[Theorem~II.10.5]{Co91}). We call $\widehat A$
the \emph{lift of $A$ relative to $\{N\}'$}. In
particular, $\widehat A$ is positive whenever $A$ is
positive.

An operator $T$ is \emph{quasinormal} if $T|T|=|T|T$,
or equivalently, $TT^*T=T^*TT$. Every subnormal
operator admits a \emph{minimal quasinormal extension}
$Q\in\ogr{\lcal}$: no proper closed $Q$-invariant
subspace of $\lcal$ containing $\hh$ carries a
quasinormal restriction of $Q$. We shall use the
following characterization (see
\cite[Theorem~1]{Cam75} and \cite[Lemma~1]{Em81}):
   \begin{align}\label{campe}
   \begin{minipage}{70ex}
\textit{If $N\in\ogr{\kk}$ is normal and
$\lcal\subseteq\kk$ is a closed $N$-invariant
subspace, then $N|_{\lcal}$ is quasinormal if and only
if $|N|(\lcal)\subseteq\lcal$.}
   \end{minipage}
   \end{align}
Finally, $T\in\ogr{\hh}$ is \emph{hyponormal} if
$\|T^*h\|\Le\|Th\|$ for every $h\in\hh$. Writing
$\mathscr{Q}$, $\mathscr{S}$, and $\mathscr{H}$ for
the classes of quasinormal, subnormal, and hyponormal
operators, respectively, we have $\mathscr{Q}\subseteq
\mathscr{S}\subseteq\mathscr{H}$
\cite[Propositions~II.1.7 and~II.4.2]{Co91}.
   \begin{dfn}\label{defq}
An operator $T\in\ogr{\hh}$ is a $B$-operator if,
relative to an orthogonal decomposition
$\hh=\hh_1\oplus\hh_2$ into nonzero subspaces, it has
the form
   \begin{align*}
T=\begin{bmatrix}V&E\\0&X\end{bmatrix},
   \end{align*}
where $V\in\ogr{\hh_1}$, $E\in\ogr{\hh_2,\hh_1}$, and
$X\in\ogr{\hh_2}$ satisfy
   \begin{gather}\label{gqb-1}
V^*V=I,
   \\\label{irxy}
V^*E=0,
   \\\label{gqb-3}
XE^*E=E^*EX \quad \big(\Leftrightarrow X|E| =
|E|X\big).
   \end{gather}
If $X$ is quasinormal \textup{(}respectively,
subnormal, hyponormal\/\textup{)}, then $T$ is said to
be of class $\mathscr{Q}$ \textup{(}respectively,
$\mathscr{S}$, $\mathscr{H}$\/\textup{)}.
   \end{dfn}
   \section{Main results}
We begin by recalling the elementary measures of a
subnormal operator and introducing their signed
counterparts. Let $S\in \ogr{\hh}$ be a subnormal
operator, and let $N\in \ogr{\kk}$ be a minimal normal
extension of $S$ with spectral measure $G_N$. The
Stone-von Neumann functional calculus (see
\cite[Theorem~5.9]{Sch12}) gives, for all integers $n
\Ge 0$ and all $f \in \hh$,
   \begin{align} \label{sdew}
\|S^n f\|^2 & = \|N^n f\|^2 = \int_{\cbb} |z|^{2n}
\langle G_N(\D z)f,f\rangle = \int_{\rbb_+}
x^n\,\D\mu_f(x),
   \end{align}
where $\mu_f$ is a finite, compactly supported Borel
measure on $\rbb_+ := [0,\infty)$ defined~by
   \begin{align} \label{gnms}
\mu_f(\varDelta) = \langle
G_N(\phi^{-1}(\varDelta))f,f\rangle, \quad \varDelta
\in \borel{\rbb_+}, \, f \in \hh,
   \end{align}
and $\phi \colon \cbb \to \rbb_+$ is defined by
$\phi(z)=|z|^2$ for $z\in\cbb$. By
\cite[Corollary~4.2]{Sch17}, the measure $\mu_f$,
called the {\em elementary measure of~$S$~at~$f$}, is
uniquely determined by $S$ and $f$ and is independent
of the choice of minimal normal extension $N$.

For an operator $T = \left[\begin{smallmatrix} V & E
\\ 0 & S \end{smallmatrix}\right] \in
\ogr{\hh_1 \oplus \hh_2}$, Boas's theorem (see
\cite[Theorem~3.11]{Sh-Ta1943}) ensures that, for
every $f\in \hh_2$, there are infinitely many finite
signed Borel measures $\rho_f$ on $\rbb_+$ such that
   \begin{align} \label{tnof}
\|T^n(0,f)\|^2 = \int_{\rbb_+} x^n \D \rho_f(x), \quad
n\Ge 0.
   \end{align}
There is at most one such measure with compact support
(see \cite[Lemma~4.1]{C-J-J-S21}). When it exists, we
call it the {\em signed elementary measure} of $T$ at
$(0,f)$; this definition does not require $T$ to be
subnormal. Exponential growth alone does not guarantee
that a real sequence admits a compactly supported
finite signed representing measure on $\rbb_+$ (see
Lemma~\ref{trkz} and Remark~\ref{negj}).

The following theorem characterizes the existence of
signed elementary measures for $B$-operators of class
$\mathscr{S}$ and expresses them in terms of the
elementary measures of $S$. It also clarifies
conditions (ii-a) and (ii-b) of Theorem~\ref{mnteo}.
   \begin{thm} \label{wkyv}
Let $T=\left[\begin{smallmatrix} V & E \\ 0 & S
\end{smallmatrix}\right] \in \ogr{\hh_1 \oplus \hh_2}$
be a $B$-operator of class $\mathscr{S}$. Then the
following conditions are equivalent{\em :}
   \begin{enumerate}
   \item[(i)] for every $f\in \hh_2$, there exists a finite
compactly supported signed Borel measure $\rho_f$ on
$\rbb_+$ satisfying \eqref{tnof},
   \item[(ii)] $\mu_{|E|f}(\{1\})=0$ for every $f\in \hh_2$,
and
   \begin{align}\label{dvcf}
\int_{\rbb_+ \setminus \{1\}}
\frac{1}{|x-1|}\,\D\mu_{|E|f}(x) < \infty, \qquad f\in
\hh_2.
   \end{align}
   \end{enumerate}
Moreover, if {\em (ii)} holds, then for all $f\in
\hh_2$ and $\varDelta \in \borel{\rbb_+}$,
   \begin{align}\label{avde}
\rho_f(\varDelta) = \int_{\rbb_+}
\frac{1}{1-x}\,\D\mu_{|E|f}(x) \, \delta_1(\varDelta)
+ \int_{\varDelta} \frac{1}{x-1}\,\D\mu_{|E|f}(x) +
\mu_f(\varDelta).
   \end{align}
   \end{thm}

We now characterize subnormality in terms of
elementary measures. The following theorem applies to
all $B$-operators of class $\mathscr{S}$.
   \begin{thm} \label{mnteo} Let $T =
\left[\begin{smallmatrix} V & E \\ 0 & S
\end{smallmatrix}\right] \in \ogr{\hh_1 \oplus \hh_2}$
be a $B$-operator of class $\mathscr{S}$. Then the
following conditions are equivalent:
   \begin{enumerate}
   \item[(i)] $T$ is subnormal,
   \item[(ii)] $T$ satisfies the following
conditions$:$
   \begin{enumerate}
   \item[(ii-a)] $\int_{\varDelta} \frac{1}{1-x}\,\D\mu_{|E|f}(x) \Le
\mu_f(\varDelta) $ for all $\varDelta \in
\borel{[0,1)}$ and all $f \in \hh_2$,
   \item[(ii-b)] $\mu_{|E|f}\big([1,\infty)\big)=0$
for all $f\in \hh_2$.
   \end{enumerate}
   \end{enumerate}
   \end{thm}
When the upper-right entry $E$ is injective,
subnormality admits a characterization in terms of a
contraction in the commutant of $|N|$, where $N$ is a
minimal normal extension of $S$.
   \begin{thm} \label{cvxw}
   Let $T = \left[\begin{smallmatrix} V & E
\\ 0 & S \end{smallmatrix}\right] \in \ogr{\hh_1 \oplus \hh_2}$
be a $B$-operator of class $\mathscr{S}$ such that $E$
is injective. Suppose that $N\in \ogr{\kk}$ is a
minimal normal extension of $S$. Then the following
conditions are equivalent{\em :}
   \begin{enumerate}
   \item[(i)]  $T$ is subnormal,
   \item[(ii)] $\|S\|\Le 1$
and there exists a contraction $R \in
\left\{|N|\right\}'$ such~that\footnote{Since, by
\cite[Corollary~II.2.17]{Co91}, $\||N|^2\|=\|N\|^2
=\|S\|^2\Le 1$, the operator $I-|N|^2$ is positive,
and hence $(I-|N|^2)^{1/2}$ is well defined.}
   \begin{align} \label{weqsa}
R (I - |N|^2)^{1/2}|_{\hh_2}= |E|.
   \end{align}
   \end{enumerate}
Moreover, if $T$ is subnormal, then $\|E\| \Le 1$,
$1\notin\sigma_{\mathrm{p}}(|N|)$ and
$G_{|N|}([1,\infty))=0$, where $G_{|N|}$ denotes the
spectral measure of $|N|$.
   \end{thm}
In the same setting, the following theorem shows that
$T$ is subnormal if and only if $T$ is a contraction.
It also provides equivalent conditions involving
Taylor spectra and spherical contractivity of pairs
associated with $(S,|E|)$ and its minimal normal and
quasinormal extensions. Section~\ref{Sec.4} and
Appendix~\ref{Sec.A} review the basic facts needed in
this article concerning the Taylor spectrum, spherical
contractions, and minimal normal and quasinormal
extensions of subnormal pairs. Recall that $(S,|E|)$
is a subnormal pair whenever $T$ is a $B$-operator of
class $\mathscr{S}$ (see Lemma~\ref{liom}). Moreover,
if $(M_1,M_2)$ is a minimal normal extension of
$(S,|E|)$ and $(Q_1,Q_2)$ is a minimal quasinormal
extension of $(S,|E|)$, then $M_2 \Ge 0$ and $Q_2 \Ge
0$. Set
   \begin{align*}
\sbd := \{(z_1,z_2)\in \cbb^2\colon |z_1|^2 + |z_2|^2
\Le 1\} \quad \text{and} \quad \sbb := \sbd \cap
\rbb_+^2.
   \end{align*}
   \begin{thm} \label{wkwsc}
Suppose that $T = \left[\begin{smallmatrix} V & E
\\ 0 & S \end{smallmatrix}\right] \in \ogr{\hh_1 \oplus \hh_2}$
is a $B$-operator of class $\mathscr{S}$ such that $E$
is injective. Then $T$ is subnormal if and only if $T$
is a contraction. Moreover, $T$ is subnormal if and
only if one of the following conditions
holds\textup{:}
   \begin{itemize}
   \item[\textup{(i)}] $\sigma(S,|E|) \subseteq
\sbd$ $($equivalently, $(S,|E|)$ is a spherical
contraction$)$,
   \item[\textup{(ii)}] $\sigma(|S|,|E|) \subseteq
\sbb$ $($equivalently, $(|S|,|E|)$ is a spherical
contraction$)$,
   \item[\textup{(iii)}] $\sigma(M_{1},M_{2}) \subseteq
\sbd$ $($equivalently, $(M_{1},M_{2})$ is a spherical
contraction$)$,
   \item[\textup{(iv)}] $\sigma(|M_{1}|,M_{2}) \subseteq
\sbb$ $($equivalently, $(|M_{1}|,M_{2})$ is a
spherical contraction$)$,
   \item[\textup{(v)}] $\sigma(Q_1, Q_2)
\subseteq \sbd$ $($equivalently, $\bigl(Q_1,
Q_2\bigr)$ is a spherical contraction$)$,
   \item[\textup{(vi)}] $\sigma(|Q_1|, Q_2)
\subseteq \sbb$ $($equivalently, $\bigl(|Q_1|,
Q_2\bigr)$ is a spherical contraction$)$,
   \end{itemize}
where $(M_1,M_2) \in \ogr{\kk}^2$ and $(Q_1,Q_2)\in
\ogr{\lcal}^2$ denote minimal normal and minimal
quasinormal extensions of $(S,|E|)$, respectively.
   \end{thm}
   To pass from the injective case to the general
case, we use the following reduction. Let $T =
\left[\begin{smallmatrix} V & E \\
0 & S \end{smallmatrix}\right] \in \ogr{\hh_1 \oplus
\hh_2}$ be a $B$-operator of class $\mathscr{S}$ with
$E \neq 0$ (the case $E=0$ is immediate, since
$T=V\oplus S$ is subnormal). Set
   \begin{align*}
\hh_{2,0}:=\overline{\ob{|E|}}, \quad E_0 :=
E|_{\hh_{2,0}}, \quad S_0:=S|_{\hh_{2,0}}, \quad
T_0:=\smat{V&E_0\\0&S_0}.
   \end{align*}
Then $\hh_{2,0}$ is invariant under $S$, $E_0 \in
\ogr{\hh_{2,0},\hh_1}$, $S_0\in \ogr{\hh_{2,0}}$ is
subnormal, and $T_0\in\ogr{\hh_1\oplus \hh_{2,0}}$ is
a $B$-operator of class $\mathscr{S}$ with $E_0$
injective. Applying Theorem~\ref{wkwsc} to $T_0$ and
using Lemma~\ref{th2z}, we obtain the following
reduction of the general case to the injective one.
   \begin{pro} \label{13posr}
Let $T = \left[\begin{smallmatrix} V & E
\\ 0 & S \end{smallmatrix}\right]$
be a $B$-operator of class $\mathscr{S}$ with $E \neq
0$. Then the following conditions are
equivalent\textup{:}
   \begin{itemize}
   \item[\textup{(i)}] $T$ is subnormal,
   \item[\textup{(ii)}] $T_0$ is subnormal,
   \item[\textup{(iii)}] $T_0$ is a contraction.
   \end{itemize}
   \end{pro}
This reduction leads to characterizations of
subnormality in terms of the Taylor spectra of pairs
associated with the second column of $T$. Set
   \begin{align*} \varXi := \sbd \cup (\cbb \times \{0\})
\quad \text{and} \quad \varXi_+ :=
\bar{\mathbb{B}}^{+}_2 \cup (\rbb_+ \times \{0\}).
   \end{align*}
   \begin{thm} \label{aswq} Let $T =
\left[\begin{smallmatrix} V & E \\ 0 & S
\end{smallmatrix}\right] \in \ogr{\hh_1 \oplus \hh_2}$
be a $B$-operator of class $\mathscr{S}$ with $E\neq
0$. Let $(M_1,M_2) \in \ogr{\kk}^2$ be a minimal
normal extension of $(S,|E|)$, and let $(Q_1,Q_2)\in
\ogr{\lcal}^2$ be a minimal quasinormal extension of
$(S,|E|)$. Then $T$ is subnormal if and only if one of
the following conditions holds\textup{:}
   \begin{itemize}
   \item[\textup{(i)}] $\sigma(S,|E|) \subseteq
\varXi$,
   \item[\textup{(ii)}] $\sigma(|S|,|E|) \subseteq
\varXi_+$,
   \item[\textup{(iii)}] $\sigma(M_{1},M_{2}) \subseteq
\varXi$,
   \item[\textup{(iv)}] $\sigma(|M_{1}|,M_{2}) \subseteq
\varXi_+$,
   \item[\textup{(v)}] $\sigma(Q_1, Q_2)
\subseteq \varXi$,
   \item[\textup{(vi)}] $\sigma(|Q_1|, Q_2)
\subseteq \varXi_+$.
   \end{itemize}
   \end{thm}
Taken together, Theorems~\ref{wkwsc} and~\ref{aswq}
and Proposition~\ref{13posr} provide a substantial
generalization of \cite[Theorem~1.2]{C-J-J-S21}.

The proofs of Theorems~\ref{wkyv} and~\ref{mnteo} are
given in Section~\ref{Sec.6}. The proofs of
Theorems~\ref{cvxw} and~\ref{wkwsc} appear in
Section~\ref{Sec.8}, while the proof of
Theorem~\ref{aswq} is given in Section~\ref{Sec.9}.
   \section{Function-theoretic preliminaries}
Let $\cbb$, $\rbb$ and $\rbb_+$ represent the sets of
all complex, real and nonnegative real numbers,
respectively. We use $\zbb$, $\zbb_+$, and $\nbb$ to
refer to the sets of integers, nonnegative integers,
and positive integers. Set
   \allowdisplaybreaks
   \begin{align*}
\tbb & =\{z\in \cbb \colon |z|=1\}, \quad \mathbb{D}
=\{z\in \cbb \colon |z|<1\}, \quad \bar\dbb =\{z \in
\cbb\colon |z| \Le 1\}.
   \end{align*}
For a Hausdorff topological space $X$, let $\borel{X}$
be the $\sigma$-algebra of all Borel subsets of $X$.
We write $\chi_{\varDelta}$ for the characteristic
function of a set $\varDelta\in\borel{X}$. The Borel
probability measure on $X$ supported on the singleton
$\{t\}$, where $t\in X$, is denoted by $\delta_t$. If
$\mu$ is a complex Borel measure on $X$, then $|\mu|$
denotes the total variation measure of $\mu$. The
closed support of a regular complex Borel measure
$\mu\colon\borel{X}\to\cbb$, or a regular positive
operator-valued Borel measure $F\colon\borel{X}\to
\ogr{\hh}$, is denoted by $\supp{\mu}$ or $\supp{F}$,
respectively. Define the linear map $\triangle\colon
\rbb^{\zbb_+}\to\rbb^{\zbb_+}$, where $\rbb^{\zbb_+}$
is the vector space of all real sequences
$\gammab=\{\gamma_n\}_{n=0}^\infty$, by
   \begin{align*}
(\triangle\gammab)_n=\gamma_{n+1}-\gamma_n,\qquad
n\in\zbb_+,\ \gammab\in\rbb^{\zbb_+}.
   \end{align*}
A sequence $\gammab = \{\gamma_n\}_{n=0}^\infty
\subseteq \rbb$ is called a {\em Stieltjes moment
sequence} if there exists a finite Borel
measure\footnote{Unless otherwise stated, measures
considered in this paper are assumed to be positive.}
$\mu$ on $\rbb_+$, referred to as a {\em representing
measure} of $\gammab$, such~that
   \begin{align*}
\gamma_n = \int_{\rbb_+} x^n \D\mu(x), \qquad n \in
\zbb_+.
   \end{align*}
   If the representing measure $\mu$ is unique, then
$\gammab$ is said to be {\em determinate}.

The {\em Stieltjes transform} of a complex Borel
measure $\rho$ on $\rbb$ is defined by
   \begin{align*}
\sfrak_{\rho} (z) = \int_{\rbb} \frac{1}{z-x} \D
\rho(x), \quad z \in \cbb \setminus \rbb.
   \end{align*}
The Stieltjes transform is a holomorphic function on
$\cbb \setminus \rbb$ (see \cite[Theorem~F.5]{Sch12}).
   \begin{lem} \label{trss}
If $\rho$ is as above and $x-1 \in L^1(|\rho|)$, then
   \begin{align*}
\sfrak_{(x-1)\, \D \rho(x)}(z)
=(z-1)\sfrak_\rho(z)-\rho(\rbb), \quad z \in \cbb
\setminus \rbb.
   \end{align*}
   \end{lem}
   \begin{proof}
Integrating the identity
   \begin{align*}
\frac{x-1}{z-x}=\frac{z-1}{z-x}-1
   \end{align*}
with respect to $\rho$ yields the desired formula.
   \end{proof}
Let $\boldsymbol{c}=\{c_n\}_{n=0}^\infty \subseteq
\cbb$ be an arbitrary sequence. By a classical result
of Boas \cite[Theorem~1]{Boa39}, there exist
infinitely many complex Borel measures $\rho$ on
$\rbb$ such~that
   \begin{align} \label{cnxn}
c_n=\int_{\rbb} x^n \D \rho(x), \quad n \in \zbb_+.
   \end{align}
Suppose now that $\rho$ in \eqref{cnxn} has compact
support. Note that such a $\rho$ is unique (see
\cite[Lemma~4.1]{C-J-J-S21}). Take $R \in (0,\infty)$
such that $\supp{\rho} \subseteq [-R, R]$. Then, for
any $z\in \cbb$ with $|z| > R$, the series
$\sum_{j=0}^\infty \left(\frac{x}{z}\right)^{j}$ is
uniformly convergent in $x$ on the interval $[-R, R]$.
Therefore, we have
    \allowdisplaybreaks
   \begin{align} \notag
\sum_{j=0}^\infty \frac{c_j}{z^{j+1}} &
\overset{\eqref{cnxn}}= \sum_{j=0}^\infty \frac{1}{z}
\int_{\rbb} \left(\frac{x}{z}\right)^{j} \D\rho(x)
   \\  \notag
& \hspace{1ex}= \frac{1}{z} \int_{\rbb}
\sum_{j=0}^\infty \left(\frac{x}{z}\right)^{j}
\D\rho(x)
   \\ \label{serxz}
& \hspace{1ex} = \int_{\rbb} \frac{1}{z-x}\,\D\rho(x),
\qquad z\in \cbb, \, |z| > R.
   \end{align}
Hence,
   \begin{align} \label{tscw}
\sum_{j=0}^\infty \frac{c_j}{z^{j+1}} =
\sfrak_{\rho}(z), \qquad z\in \cbb\setminus \rbb,\ |z|
> R.
   \end{align}
Note also that
   \begin{align*}
\sum_{j=0}^\infty \frac{c_j}{z^{j+1}} =
\varphi\left(\frac{R}{z}\right), \qquad z\in \cbb,\
|z|
> R,
   \end{align*}
where $\varphi$ is a holomorphic function on $\dbb$
defined by the power series
   \begin{align*}
\varphi(w)=\sum_{j=0}^\infty \frac{c_j}{R^{j+1}}
w^{j+1}, \quad w \in \dbb.
   \end{align*}
   \section{\label{Sec.4}Operator-theoretic preliminaries}
We begin with the following general observation.
   \begin{lem} \label{comx}
Let $A,B,R\in\ogr{\kk}$, where $A$ and $B$ are
selfadjoint and $B$ is injective. Suppose $AB=BA$ and
$RB=A$. Then $R$ is selfadjoint and $RB=BR$.
   \end{lem}
   \begin{proof}
Since $B$ is selfadjoint and injective, we obtain
$\overline{\ob{B}}=\kk$. The identity $AB=BA$ implies
$AB^{-1}\subseteq B^{-1}A$. Moreover, $RB=A=A^*=BR^*$,
so $\ob{A}\subseteq\ob{B}$ and $B^{-1}A=R^*$. Thus
   \begin{align*}
R|_{\ob{B}}=AB^{-1}\subseteq B^{-1}A=R^*,
   \end{align*}
which implies that $R=R^*$ and hence $RB=BR$.
   \end{proof}
We now turn to some necessary concepts from
multivariate operator theory. Given
$\tbold=(T_1,T_2)\in \ogr{\hh}^2$, we write $\tbold^*=
(T_1^*,T_2^*)$ and
$\tbold^{\alphab}=T_1^{\alpha_1}T_2^{\alpha_2}$ for
all $\alphab=(\alpha_1,\alpha_2)\in \zbb_+^2$. The
pair $\tbold$ is called a {\em spherical contraction}
if $T_1^*T_1 + T_2^*T_2 \Le I$. Clearly, $\tbold$ is a
spherical contraction if and only if $(|T_1|,|T_2|)$
is a spherical contraction, and if $\tbold$ is a
spherical contraction, then both components $T_1$ and
$T_2$ are contractions. Following \cite{Ath88}, we say
that $\tbold$ is {\em hyponormal} if
   \begin{align*}
[\tbold^*,\tbold]:=\begin{bmatrix} T_1^*T_1 - T_1T_1^*
& T_2^*T_1 - T_1T_2^*
   \\[1ex]
T_1^*T_2 - T_2T_1^* & T_2^*T_2 - T_2T_2^*
   \end{bmatrix} \Ge 0.
   \end{align*}
It follows from the definition that if $\mathbf{T}$ is
hyponormal, then both $T_1$ and $T_2$ are hyponormal.

Suppose that $\tbold=(T_1,T_2)\in \ogr{\hh}^2$ is a
pair of commuting operators. Denote by
$\sigma(\tbold)=\sigma(T_1,T_2)$ the Taylor spectrum
of $\tbold$ (see \cite{Tay70}). Using
\cite[Theorem~1.1]{Vas78}, we deduce that if
$(S_1,S_2)\in \ogr{\kk}^2$ is another pair of
commuting operators, then
   \begin{align} \label{ossp}
\sigma\big((T_1, T_2) \oplus (S_1, S_2)\big) =
\sigma(T_1, T_2) \cup \sigma(S_1, S_2),
   \end{align}
where the pair $(T_1, T_2) \oplus (S_1, S_2) \in
\ogr{\hh \oplus \kk}^2$ of commuting operators is
defined~by
   \begin{align*}
(T_1, T_2) \oplus (S_1, S_2) := (T_1 \oplus S_1, T_2
\oplus S_2).
   \end{align*}
The pair $\tbold$ is said to be {\em subnormal} if
there exists a Hilbert space $\kk$ and a pair
$\nbold=(N_1,N_2)\in \ogr{\kk}^2$ of commuting normal
operators such that $\hh \subseteq \kk$ and $T_j
\subseteq N_j$ for $j=1,2$; such a pair $\nbold$ is
called a {\em normal extension} of $\tbold$. If $\kk$
is the only closed subspace of $\kk$ containing $\hh$
and reducing $N_j$ for $j=1,2$, then $\nbold$ is
called a {\em minimal normal extension} of $\tbold$.
We recall the following for later use.
   \begin{align} \label{mnex}
   \begin{minipage}{70ex}
\textit{A normal extension $\nbold$ of $\tbold$ is
minimal $\iff$ $\kk=\bigvee_{\alpha\in \zbb_+^2}
\nbold^{*\alpha} \hh$.}
   \end{minipage}
   \end{align}
If $\nbold=(N_1,N_2)\in \ogr{\kk}^2$ and
$\mbold=(M_1,M_2)\in \ogr{\lcal}^2$ are minimal normal
extensions of a subnormal pair $\tbold$, then $\nbold$
and $\mbold$ are {\em $\hh$-unitarily equivalent},
that is, there exists a unitary operator $U\in
\ogr{\kk,\lcal}$ such that $U|_{\hh}=I_{\hh}$ and
$UN_j=M_jU$ for $j=1,2.$ Recall the spectral inclusion
theorem for subnormal pairs (see
\cite[Theorem~1]{Cur81}; see also
\cite[Theorem]{Put84}):
   \begin{align} \label{pytr}
   \begin{minipage}{65ex}
{\em If $\nbold$ is a minimal normal extension of
$\tbold$, then $\sigma(\nbold) \subseteq
\sigma(\tbold)$.}
   \end{minipage}
   \end{align}
Subnormal pairs are known to be hyponormal (see
\cite[Proposition~2]{Ath88}).

The following two lemmas play a crucial role in this
paper.
   \begin{lem} \label{liom}
Let $S\in \ogr{\hh}$ be a subnormal operator that
commutes with a positive operator $A\in \ogr{\hh}$,
let $N\in \ogr{\kk}$ be a minimal normal extension of
$S$, and let $\widehat{A}$ be the lift of $A$ relative
to $\{N\}^{\prime}$. Then the pair $(S,A)$ is
subnormal, and $(N,\widehat{A})$ is a minimal normal
extension of $(S,A)$.
   \end{lem}
   \begin{proof}
Clearly, $(N,\widehat{A})$ is a pair of commuting
normal operators extending $(S,A)$; hence $(S,A)$ is
subnormal. Since $N$ is a minimal normal extension of
$S$, it follows from \eqref{mnex} that
$(N,\widehat{A})$ is minimal.
   \end{proof}
We now turn to the case where the operator $S$ in
Lemma~\ref{liom} is hyponormal.
   \begin{lem}[{\cite[Remark~4]{Ath88}}]  \label{lohp}
Let $S\in \ogr{\hh}$ be a hyponormal operator that
commutes with a positive operator $A\in \ogr{\hh}$.
Then the pair $(S,A)$ is hyponormal.
   \end{lem}
We next establish a result that relates the
contractivity of $B$-operators to both the spherical
and spectral contractivity of the associated pair
$(S,|E|)$, and to the spectral contractivity of a
minimal normal extension $\nbold=(N_1, N_2)$ of
$(S,|E|)$. Recall that for a pair $\tbold=(T_1,T_2)\in
\ogr{\hh}^2$ of commuting operators, the spectral
radius $r(\tbold)$ is defined by
   \begin{align*}
r(\tbold)=r(T_1,T_2) = \max \Big\{\sqrt{|z_1|^2 +
|z_2|^2} \colon (z_1, z_2) \in \sigma(\tbold)\Big\},
   \end{align*}
and the elementary operator $M_{\tbold} \colon
\ogr{\hh} \to \ogr{\hh}$ is defined by
   \begin{align} \label{mtvs}
M_{\tbold}(X) = T_1^*XT_1 + T_2^*XT_2, \quad X\in
\ogr{\hh}.
   \end{align}
   \begin{lem} \label{hcps}
Let $\tbold=(T_1,T_2)\in \ogr{\hh}^2$ be a hyponormal
pair of commuting operators such that $|T_1|$ and
$|T_2|$ commute. Then the following conditions are
equivalent{\em :}
   \begin{enumerate}
   \item[(i)] $\sigma(\tbold) \subseteq
\sbd$,
   \item[(ii)] $\sigma(|T_1|,|T_2|) \subseteq
\sbb$,
   \item[(iii)] $\tbold$ is a spherical contraction.
   \end{enumerate}
   \end{lem}
   \begin{proof}
It follows from \cite[Lemma~3.10]{Chav-Sh17} (see also
\eqref{mtvs}) that
   \begin{gather*}
r(T_1,T_2) = \|M_{(T_1,T_2)}(I)\|^{1/2} = \|T_1^*T_1 +
T_2^*T_2\|^{1/2},
   \\
r(|T_1|,|T_2|) = \|M_{(|T_1|,|T_2|)}(I)\|^{1/2} =
\|T_1^*T_1 + T_2^*T_2\|^{1/2}.
   \end{gather*}
This implies that conditions (i)-(iii) are equivalent.
   \end{proof}
Note that the implication (iii)$\Rightarrow$(i) holds
without assuming the hyponormality of $\tbold$ (see
\cite[Lemma~2.1]{Ch-15}) and that the equivalence
(ii)$\Leftrightarrow$(iii) can be proved by applying
the spectral theorem for a pair of bounded commuting
selfadjoint operators.
   \begin{lem} \label{cxzq}
Let $S\in\ogr\hh$ be a subnormal operator that
commutes with a positive injective operator
$A\in\ogr\hh$. Then the following conditions are
equivalent\textup{:}
   \begin{enumerate}
   \item[(i)] $\sigma(S,A)\subseteq \varXi$,
   \item[(ii)] $\sigma(S,A)\subseteq\sbd$.
   \end{enumerate}
If, in addition, $S$ is positive, then conditions
\textup{(i)}-\textup{(ii)} remain equivalent when
$\varXi$ and $\sbd$ are replaced by $\varXi_+$ and
$\sbb$, respectively.
   \end{lem}
   \begin{proof}
Let $(N,\widehat A)$ be as in Lemma~\ref{liom}. By
Lemma~\ref{lohp}, the pairs $(S, A)$ and $(N,\widehat
A)$ are hyponormal. Denote by $G_N$ and $G_{\widehat
A}$ the spectral measures of $N$ and $\widehat A$,
respectively, and by $G$ the joint spectral measure of
$(N,\widehat A)$. Recall that $G$ is the product of
the measures $G_N$ and $G_{\widehat A}$ (see
\cite[Subsection~5.5.1]{Sch12}). Since $A$ is
injective, so is $\widehat A$ (see
\cite[Lemma~3.2]{C-J-J-S24}). Thus, $G_{\widehat
A}(\{0\})=0$ (see \cite[Theorem~6.6.2]{Bi-So87}),
which implies that
   \begin{equation}\label{ilkp}
G(\cbb\times\{0\})=G_N(\cbb)G_{\widehat A}(\{0\})=0.
   \end{equation}
It follows from \cite[Theorem~2.1(i)]{C-J-J-S21} (see
also Footnote~3 therein) that
   \begin{align}\label{spgs}
\supp G=\sigma(N,\widehat A).
   \end{align}

\textup{(i)}$\Rightarrow$\textup{(ii)} Suppose that
$\sigma(S,A)\subseteq\varXi$. The spectral inclusion
\eqref{pytr} yields
   \begin{align}\label{sass}
\sigma(N,\widehat A)\subseteq\sigma(S,A)\subseteq
\varXi.
   \end{align}
Since $\varXi\setminus\sbd\subseteq\cbb\times\{0\}$,
we deduce from \eqref{ilkp}, \eqref{spgs}, and
\eqref{sass} that
   \begin{align*}
G(\cbb^2\setminus\sbd)
&=G(\cbb^2\setminus\varXi)+G(\varXi\setminus\sbd)\\
&\Le G(\cbb^2\setminus\supp G) +G(\cbb\times\{0\})=0.
   \end{align*}
Since $\sbd$ is closed, we have
   \begin{align*}
\sigma(N,\widehat A)=\supp G\subseteq\sbd.
   \end{align*}
By Lemma~\ref{hcps} (or directly by the spectral
theorem), $N^*N+\widehat A^{\,2}\Le I$. Restricting
this quadratic-form inequality to $\hh$ gives
$S^*S+A^2\Le I$. Applying Lemma~\ref{hcps} again, we
conclude that $\sigma(S,A)\subseteq\sbd$.

\textup{(ii)}$\Rightarrow$\textup{(i)} This is
immediate.

Finally, if $S$ is positive, the assertion follows
from the preceding argument and the projection
property of the Taylor spectrum (see \cite[p.\
21]{Cur88}).
   \end{proof}
   \begin{lem}  \label{snqm}
Let $S\in \ogr{\hh}$ be a subnormal operator that
commutes with a positive operator $A\in \ogr{\hh}$,
let $N\in \ogr{\kk}$ be a minimal normal extension of
$S$, and let $\widehat{A}$ be the lift of $A$ relative
to $\{N\}^{\prime}$. Then $\lcal := \bigvee_{n\Ge 0}
|N|^n \hh$ is invariant under $N$ and $\widehat{A}$,
$(N|_{\lcal}, \widehat{A}|_{\lcal})$ is a minimal
quasinormal extension of $(S, A)$, and the following
conditions are equivalent\textup{:}
   \begin{enumerate}
   \item[(i)] $\sigma(S,A)\subseteq \sbd$,
   \item[(ii)] $\sigma(N|_{\lcal}, \widehat{A}|_{\lcal})
\subseteq\sbd$,
   \item[(iii)] $\sigma(N, \widehat{A})
\subseteq\sbd$.
   \end{enumerate}
   \end{lem}
   \begin{proof}
By Lemma~\ref{liom}, $\nbold:= (N, \widehat{A})$ is a
minimal normal extension of $\sbold:=(S, A)$. Since
$\lcal = \bigvee_{m,n \Ge 0} |N|^m \widehat{A}^n \hh$,
it follows from Appendix~\ref{Subs.11.2} that $\lcal$
is invariant under $N$ and $\widehat{A}$,
$\qbold:=(N|_{\lcal}, \widehat{A}|_{\lcal})$ is a
minimal quasinormal extension of $\sbold$, and
$\nbold$ is a minimal normal extension of $\qbold$.
Moreover, $N$ is a minimal normal extension of
$N|_{\lcal}$ and $\widehat{A}$ is the lift of
$\widehat{A}|_{\lcal}$ relative to $\{N\}^{\prime}$.

The implication (i)$\Rightarrow$(iii) follows from
Lemma~\ref{liom} and \eqref{pytr}. The converse
implication (iii)$\Rightarrow$(i) follows from
\cite[Theorem~2]{Cur81} and $\mathrm{co}(\sbd)=\sbd$,
because
   \begin{align*}
\sigma(S,A) \subseteq \mathrm{co}(\sigma(N,
\widehat{A})) \overset{\mathrm{(iii)}}{\subseteq}
\sbd,
   \end{align*}
where $\mathrm{co}(X)$ denotes the polynomially convex
hull of $X \subseteq \cbb^2$.

Applying the equivalence (i)$\Leftrightarrow$(iii) to
$(\qbold, \nbold)$ in place of $(\sbold, \nbold)$, we
obtain the equivalence (ii)$\Leftrightarrow$(iii).
   \end{proof}
   \section{Subnormality via moment techniques}
This section contains three lemmas that prepare the
ground for the further study of the subnormality of
$B$-operators.
   \begin{lem}
   Let $T =
\left[\begin{smallmatrix} V & E \\ 0 & S
\end{smallmatrix}\right] \in \ogr{\hh_1 \oplus \hh_2}$
be a $B$-operator. Then
   \begin{align} \label{tstn}
T^{*n} T^n = \begin{bmatrix} I & 0 \\ 0 &
\varOmega_{n}
\end{bmatrix},  \quad n \in \zbb_+,
   \end{align}
where
   \begin{align} \label{vyom}
\varOmega_{n} = \varOmega_{n}(T) :=
\begin{cases}
I & \text{if } n = 0, \\
E^*E \left( \sum_{j=0}^{n-1} S^{*j} S^j \right) +
S^{*n} S^n & \text{if } n \in \nbb.
\end{cases}
   \end{align}
In particular, the following identity holds{\em :}
   \begin{equation} \label{eq:Tn-norm}
\|T^n(0,f)\|^2=\is{\varOmega_{n}
f}{f}=\sum_{j=0}^{n-1}\|S^j |E| f\|^2 + \|S^n f\|^2,
\quad n \in \nbb,\, f \in \hh_2.
   \end{equation}
   \end{lem}
   \begin{proof}
To prove \eqref{tstn}, we proceed by induction on $n$.
The case $n=0$ is trivial. Suppose that the formula
holds for some $n\Ge 0$. Then, by \eqref{gqb-1}, we
have
   \allowdisplaybreaks
   \begin{align*}
T^{* (n+1)} T^{n+1} &= T^* (T^{*n} T^n)  T  \\
&=
\begin{bmatrix}
V^* & 0 \\
E^* & S^*
\end{bmatrix}
\begin{bmatrix}
I & 0 \\
0 & \varOmega_n
\end{bmatrix}
\begin{bmatrix}
V & E \\
0 & S
\end{bmatrix} \\
&=
\begin{bmatrix}
V^* & 0 \\
E^* & S^*\varOmega_{n}
\end{bmatrix}
\begin{bmatrix}
V & E \\
0 & S
\end{bmatrix} \\
&=
\begin{bmatrix}
I & 0 \\
0 & E^* E + S^* \varOmega_{n} S
\end{bmatrix}.
\end{align*}
This implies that
$$
\varOmega_{n+1} = E^* E + S^* \varOmega_{n} S,
$$
and thus by the induction hypothesis and
\eqref{gqb-3}, we get
   \allowdisplaybreaks
   \begin{align*}
\varOmega_{n+1} &= E^* E + S^* \bigg(E^*E
\bigg(\sum_{j=0}^{n-1}
S^{*j} S^j \bigg) + S^{*n} S^n \bigg) S \\
&= E^* E + E^*E \bigg(\sum_{j=1}^{n} S^{*j} S^j \bigg)
+ S^{*(n+1)} S^{n+1} \\
&= E^* E \bigg( \sum_{j=0}^{n} S^{*j} S^j \bigg) +
S^{*(n+1)} S^{n+1}.
   \end{align*}
The ``in particular'' part follows from the above
using \eqref{gqb-3}.
   \end{proof}
The following lemma gives an integral representation
of the generating function of
$\{\|T^n(0,f)\|^2\}_{n=0}^\infty$ via the elementary
measures of $S$.
   \begin{lem}\label{lem:gen}
Let $T = \left[\begin{smallmatrix} V & E
\\ 0 & S \end{smallmatrix}\right] \in \ogr{\hh_1 \oplus \hh_2}$
be a $B$-operator of class $\mathscr{S}$. Fix
$f\in\hh_2$. Set $g:=|E|f$ and $c_n:=\|T^n(0,f)\|^2$
for $n\in \zbb_+$. Then
   \begin{equation}\label{eq:gen}
\sum_{n=0}^\infty \frac{c_n}{z^{n+1}} = \int_{\rbb_+}
\frac{1}{z-x}\,\D\mu_f(x) + \frac{1}{z-1}\int_{\rbb_+}
\frac{1}{z-x}\,\D\mu_g(x),
   \end{equation}
for every $z \in \cbb$ such that $|z| > R$, where $R
:= \max \bigl \{1, \sup \supp{\mu_f}, \sup
\supp{\mu_g} \bigr\}$, and $\mu_h$ denotes the
elementary measure of $S$ at $h\in \hh_2$ $($see
\eqref{sdew}$)$.
   \end{lem}
   \begin{proof}
It follows from \eqref{eq:Tn-norm} that
   \begin{align}  \label{cnjx}
c_n= \|S^n f\|^2+\sum_{j=0}^{n-1}\|S^j g\|^2, \quad
n\in\nbb.
   \end{align}
Clearly, $c_0=\|f\|^2$. Hence, for any real $t > R$,
we have
\[
\sum_{n=0}^\infty \frac{c_n}{t^{n+1}} =
\sum_{n=0}^\infty \frac{\|S^n f\|^2}{t^{n+1}} +
\sum_{n=1}^\infty
\frac{1}{t^{n+1}}\sum_{j=0}^{n-1}\|S^j g\|^2.
\]
Using Tonelli's theorem for nonnegative terms, we
obtain
   \allowdisplaybreaks
   \begin{align} \notag
\sum_{n=1}^\infty
\frac{1}{t^{n+1}}\sum_{j=0}^{n-1}\|S^j g\|^2 &
=\sum_{j=0}^\infty \|S^j g\|^2\sum_{n=j+1}^\infty
\frac{1}{t^{n+1}}
   \\ \notag
&=\sum_{j=0}^\infty \|S^j g\|^2\cdot
\frac{1}{t^{j+2}}\cdot \frac{1}{1-\frac1t}
   \\ \notag
& =\frac{1}{t-1}\sum_{j=0}^\infty \frac{\|S^j
g\|^2}{t^{j+1}}
   \\ \label{toli}
& =\frac{1}{t-1}\int_{\rbb_+} \frac{1}{t-x}\,\D
\mu_g(x) < \infty, \qquad t>R,
   \end{align}
where the last equality follows from \eqref{serxz}
with $\rho=\mu_g$. This shows that the double series
in \eqref{toli} is absolutely convergent.
Consequently, assuming $|z|>R$, we may repeat the
above argument with $z$ in place of $t$, then apply
\eqref{serxz} with $\rho=\mu_f$, and finally combine
the result with \eqref{cnjx} to obtain \eqref{eq:gen}.
   \end{proof}
Moment problems are fundamentally connected to
subnormality, as highlighted by Lambert's celebrated
theorem.
   \begin{thm}[{\cite{Lam76,St-SzII}}]
\label{Lams} An operator $T\in \ogr{\hh}$ is subnormal
if and only if the sequence
$\{\|T^nh\|^2\}_{n=0}^{\infty}$ is a Stieltjes moment
sequence for every $h\in \hh$.
   \end{thm}
This connection enables us to characterize the
subnormality of $B$-operators in terms of Stieltjes
moment sequences.
   \begin{lem} \label{vyox}
Let $T = \left[\begin{smallmatrix} V & E \\ 0 & S
\end{smallmatrix}\right] \in \ogr{\hh_1 \oplus \hh_2}$
be a $B$-operator. Then $T$ is subnormal if and only
if the sequence $\left\{ \langle \varOmega_n f, f
\rangle \right\}_{n=0}^\infty $ is a Stieltjes moment
sequence for every $f \in \hh_2$, where $\varOmega_n$
is defined as in \eqref{vyom}.
   \end{lem}
   \begin{proof}
It follows from \eqref{tstn} that for $f = f_1 \oplus
f_2 \in \hh_1 \oplus \hh_2$, we have
   \begin{align*}
\langle T^{*n} T^n f, f \rangle = \| f_1 \|^2 +
\langle \varOmega_n f_2, f_2 \rangle, \quad n \in
\zbb_+.
   \end{align*}
This implies that the sequence $ \left\{ \langle
T^{*n} T^n f, f \rangle \right\}_{n=0}^\infty $ is a
Stieltjes moment sequence for every $f \in \hh$ if and
only if $ \left\{ \langle \varOmega_n f_2, f_2 \rangle
\right\}_{n=0}^\infty $ is a Stieltjes moment sequence
for every $f_2 \in \hh_2$. An application of
Theorem~\ref{Lams} completes the proof.
   \end{proof}
   \section{\label{Sec.6}Subnormality via elementary measures}
   \begin{proof}[Proof of Theorem~\ref{wkyv}]
   Fix $f\in \hh_2$, and write $c^f_n=\|T^n(0,f)\|^2$
for all $n\in \zbb_+$. Since $S$ is subnormal, it
follows from \eqref{sdew} and \eqref{eq:Tn-norm} that
   \allowdisplaybreaks
   \begin{align} \label{msaz}
c^f_n = \sum_{j=0}^{n-1}
\int_{\rbb_+}x^j\,\D\mu_{|E|f}(x) + \int_{\rbb_+}x^n\,
\D\mu_f(x), \quad n \in \nbb,
   \end{align}
and $c^f_0 = \int_{\rbb_+} x^0 \, \D\mu_f(x)$. Both
$\mu_f$ and $\mu_{|E|f}$ are finite compactly
supported positive Borel measures.

(i)$\Rightarrow$(ii) Suppose that $\rho_f$ is a finite
compactly supported signed Borel measure on $\rbb_+$
such that
   \begin{align} \label{cfni}
c^f_n = \int_{\rbb_+} x^n \D \rho_f(x), \quad n\in
\zbb_+.
   \end{align}
Taking first differences, we obtain
   \begin{align*}
\int_{\rbb_+} x^n (x-1)\,& \D\rho_f(x)
=c^f_{n+1}-c^f_n
   \\
&\hspace{-1ex}\overset{\eqref{msaz}}= \int_{\rbb_+}x^n
(x-1)\,\D\mu_f(x) +\int_{\rbb_+} x^n \,
\D\mu_{|E|f}(x), \quad n \in \zbb_+.
   \end{align*}
Therefore, $(x-1)\, \D (\rho_f-\mu_f)(x)$ and
$\mu_{|E|f}$ are finite compactly supported signed
Borel measures on $\rbb_+$ with the same moments.
Hence, by \cite[Lemma~4.1]{C-J-J-S21}, we~have
   \begin{align} \label{xronu}
\mu_{|E|f}(\varDelta) = \int_{\varDelta} (x-1)\, \D
(\rho_f-\mu_f)(x), \quad \varDelta \in \borel{\rbb_+}.
   \end{align}
First, evaluating the above identity on $\{1\}$ shows
that
   \begin{align*}
\mu_{|E|f}(\{1\})=0.
   \end{align*}
Using \cite[Theorems~6.12 and 6.13]{Rud87}, we deduce
from \eqref{xronu} that
   \begin{align*}
\mu_{|E|f}(\varDelta) = \int_{\varDelta} |x-1| \, \D
|\rho_f-\mu_f|(x), \quad \varDelta \in \borel{\rbb_+}.
   \end{align*}
This implies that
   \begin{align*}
\int_{\rbb_+\setminus \{1\}} \frac{1}{|x-1|} \, \D
\mu_{|E|f}(x) = |\rho_f-\mu_f|(\rbb_+\setminus \{1\})
< \infty.
   \end{align*}
This, together with \eqref{xronu}, yields
   \begin{align*}
\int_{\varDelta} \frac{1}{x-1} \, \D \mu_{|E|f}(x) =
\rho_f(\varDelta\setminus
\{1\})-\mu_f(\varDelta\setminus \{1\}), \quad
\varDelta \in \borel{\rbb_+}.
   \end{align*}
Consequently, for every $\varDelta \in
\borel{\rbb_+}$,
   \begin{align} \notag
\rho_f(\varDelta) & = \rho_f(\{1\}) \delta_1
(\varDelta) + \rho_f(\varDelta\setminus \{1\})
   \\ \label{dwqp}
& = (\rho_f-\mu_f)(\{1\}) \, \delta_1 (\varDelta) +
\mu_f(\varDelta) + \int_{\varDelta} \frac{1}{x-1} \,
\D \mu_{|E|f}(x).
   \end{align}
Taking $\varDelta=\rbb_+$ and using the equalities
   \begin{align} \label{trrci}
\mu_f(\rbb_+)\overset{\eqref{sdew}}=\|f\|^2
\overset{\eqref{tnof}}=\rho_f(\rbb_+),
   \end{align}
we obtain
   \begin{align*}
(\rho_f-\mu_f)(\{1\}) = \int_{\rbb_+} \frac{1}{1-x} \,
\D \mu_{|E|f}(x).
   \end{align*}
This, together with \eqref{dwqp}, implies
\eqref{avde}, which proves the ``moreover'' part.

\textup{(ii)}$\Rightarrow$\textup{(i)} Suppose that
\textup{(ii)} holds. Then \eqref{avde} defines a
finite compactly supported signed Borel measure
$\rho_f$ on $\rbb_+$. Moreover,
   \begin{align*}
\rho_f(\rbb_+)=\mu_f(\rbb_+)=\|f\|^2=c^f_0,
   \end{align*}
and
   \allowdisplaybreaks
   \begin{align*}
\int_{\rbb_+}x^n\,\D\rho_f(x) &=
\int_{\rbb_+}x^n\,\D\mu_f(x) +
\int_{\rbb_+\setminus\{1\}}
\frac{x^n-1}{x-1}\,\D\mu_{|E|f}(x)
   \\
&= \int_{\rbb_+}x^n\,\D\mu_f(x) + \sum_{j=0}^{n-1}
\int_{\rbb_+}x^j\,\D\mu_{|E|f}(x)
   \\
& \hspace{-1ex} \overset{\eqref{msaz}}= c^f_n, \qquad
n\in\nbb,
   \end{align*}
Thus \eqref{cfni} holds. This completes the proof.
   \end{proof}
The alternative proof of (i)$\Rightarrow$(ii) in
Theorem~\ref{wkyv} given below uses the Stieltjes
transform to establish the key identity \eqref{xronu};
the remaining steps are as in the preceding proof.
This approach also yields further information about
the Stieltjes transform of $\rho_f-\mu_f$ (see
\eqref{s2sm}).
   \begin{proof}[Proof II of
Theorem~\ref{wkyv}] We sketch the proof of
(i)$\Rightarrow$(ii). Fix $f\in\hh_2$ and set
$g=|E|f$. Using Lemma~\ref{lem:gen} and applying
\eqref{tscw} with $\rho=\rho_f$, we obtain
   \begin{align*}
\sfrak_{\rho_f}(z) = \sfrak_{\mu_f}(z) + \frac{1}{z-1}
\sfrak_{\mu_g}(z), \qquad z\in \cbb\setminus \rbb,\
|z| > R,
   \end{align*}
where $R\in [1,\infty)$ is chosen so that
$\supp{\mu_f} \cup \supp{\mu_g} \cup \supp{\rho_f}
\subseteq [-R, R]$. The identity theorem for
holomorphic functions then yields
   \begin{align} \label{s2sm}
(z-1) \sfrak_{\sigma_f}(z) = \sfrak_{\mu_g}(z), \quad
z\in \cbb\setminus \rbb,
   \end{align}
where $\sigma_f:=\rho_f-\mu_f$ is a finite signed
Borel measure on $\rbb_+$ with compact support.
Applying \cite[F.5, p.\ 412]{Sch12} to the measures
$\mu_g$ and $\sigma_f$, we obtain
   \begin{align*}
\mu_g(\{1\}) = \lim_{\varepsilon \to 0+} \I
\varepsilon \, \sfrak_{\mu_g} (1+\I \varepsilon)
\overset{\eqref{s2sm}}= \lim_{\varepsilon \to 0+} (\I
\varepsilon)^2 \, \sfrak_{\sigma_f} (1+\I
\varepsilon)=0.
   \end{align*}
Next, Lemma~\ref{trss} gives
   \begin{align*}
\sfrak_{(x-1)\, \D \sigma_f(x)}(z) & = (z-1)
\sfrak_{\sigma_f}(z) - (\rho_f(\rbb_+)-\mu_f(\rbb_+))
   \\
& \hspace{-1ex} \overset{\eqref{trrci}} =(z-1)
\sfrak_{\sigma_f}(z)
   \\
& \hspace{-1ex} \overset{\eqref{s2sm}} =
\sfrak_{\mu_g}(z), \quad z\in \cbb\setminus \rbb.
   \end{align*}
It follows from the Stieltjes-Perron inversion formula
(see \cite[Theorem~A.13]{Sch17}) that \eqref{xronu}
holds.
   \end{proof}
   \begin{proof}[Proof of Theorem~\ref{mnteo}]
\textup{(i)}$\Rightarrow$\textup{(ii)} Suppose that
$T$ is subnormal. It follows from \eqref{eq:Tn-norm}
and Lemma~\ref{vyox} that the condition \textup{(i)}
of Theorem~\ref{wkyv} holds, and that all the measures
$\rho_f$ are positive. By this theorem, \eqref{dvcf}
and \eqref{avde} hold, and
   \begin{align} \label{mmef}
\mu_{|E|f}(\{1\})=0, \quad f \in \hh_2.
   \end{align}
Since $\rho_f\big|_{\borel{[0,1)}} \Ge 0$ and
$\rho_f(\{1\}) \Ge 0$ for all $f\in\hh_2$, we deduce
from \eqref{avde} that \textup{(ii-a)} and
\textup{(ii-b$^\prime$)} are valid, respectively,
where
   \begin{align*}
\textup{(ii-b$^{\prime}$)} \quad & \int_{(1,\infty)}
\frac{1}{x-1}\, \D\mu_{|E|f}(x) \Le \int_{[0,1)}
\frac{1}{1-x}\, \D\mu_{|E|f}(x) + \mu_f(\{1\}) \text{
for all } f \in \hh_2.
   \end{align*}
Now observe that, by \eqref{sdew}, we have
   \begin{align*}
\int_{\mathbb{R}_+} x^k \, \D \mu_{S^n f}(x) =
\|S^{k+n} f\|^2 = \int_{\rbb_+} x^k \, x^n \D
\mu_{f}(x), \qquad k\in \zbb_+, \, n\in \zbb_+.
   \end{align*}
Since the measures $\D\mu_{S^k f}(x)$ and $x^k \D
\mu_{f}(x)$ are finite compactly supported, the
corresponding Stieltjes moment problem is determinate,
and hence
   \begin{align} \label{dmxk}
\D\mu_{S^k f}(x)=x^k\, \D\mu_f(x), \qquad f\in \hh_2.
   \end{align}
It follows from Definition~\ref{defq} that $S^n
|E|=|E|S^n$ for all $n\in \zbb_+$, which implies~that
   \allowdisplaybreaks
   \begin{align*}
\int_{(1,\infty)} \frac{x^n}{x-1}\, \D\mu_{|E| f}(x) &
\overset{\eqref{dmxk}}= \int_{(1,\infty)}
\frac{1}{x-1}\, \D\mu_{S^n|E| f}(x)
   \\
& \hspace{1ex}= \int_{(1,\infty)} \frac{1}{x-1}\,
\D\mu_{|E| S^nf}(x)
   \\
& \hspace{-.55ex}
\overset{\text{(ii-b$^{\prime}$)}}{\Le} \int_{[0,1)}
\frac{1}{1-x}\, \D\mu_{S^n |E| f}(x) + \mu_{S^n
f}(\{1\})
   \\
&\overset{\mathrm{\eqref{dmxk}}}= \int_{[0,1)}
\frac{x^n}{1-x}\, \D\mu_{|E|f}(x) + \mu_{f}(\{1\}),
   \end{align*}
for all $n\in \zbb_+$ and $f \in \hh_2$. But
\[
\frac{x^n}{x-1}\nearrow +\infty \quad \text{as }
n\to\infty \quad \text{for all } x\in(1,\infty),
\]
and
\[
\frac{x^n}{1-x}\to 0 \quad \text{as } n\to\infty \quad
\text{for all } x\in[0,1).
\]
Moreover (cf.\ (ii-a)),
\[
0\Le \frac{x^n}{1-x} \Le \frac{1}{1-x} \text{ for all
} x \in [0,1), \text{ and } \int_{[0,1)} \frac{1}{1-x}
\, \D \mu_{|E| f}(x) < \infty.
\]
Therefore, applying the Lebesgue monotone convergence
theorem and the dominated convergence theorem, we
obtain $\mu_{|E|f}\big((1,\infty)\big)=0$ for all
$f\in \hh_2$. This, together with \eqref{mmef}, gives
(ii-b).

\textup{(ii)}$\Rightarrow$\textup{(i)} Conditions
(ii-a) and (ii-b) imply the condition (ii) of Theorem~
\ref{wkyv}. By this theorem, its condition (i) is
valid, and the measures $\rho_f$ are given by
\eqref{avde}. By \eqref{avde},
$\rho_f\big|_{\borel{(1,\infty)}} \Ge 0$ for all
$f\in\hh_2$. In turn, applying (ii-a) and (ii-b)
again, we deduce from \eqref{avde} that
$\rho_f\big|_{\borel{[0,1)}} \Ge 0$ and $\rho_f(\{1\})
\Ge 0$ for all $f\in\hh_2$. Thus each $\rho_f$ is
positive. Hence, by \eqref{tnof}, \eqref{eq:Tn-norm},
and Lemma~\ref{vyox}, $T$ is subnormal.
   \end{proof}
In connection with Boas's theorem (cf.\
\cite[Theorem~1]{Boa39}), we show that there exist
real sequences of exponential growth that cannot be
represented by a compactly supported finite signed
Borel measure on $\rbb_+$ (see Remark~\ref{negj}).
   \begin{lem} \label{trkz}
Let $m \Ge 1$ be an integer, and let
$\boldsymbol{c}=\{c_n\}_{n=0}^\infty \subseteq \mathbb
R$ be a nonconstant sequence such that $\triangle^m
\boldsymbol{c} \equiv 0$. Then there exists no finite
signed Borel measure $\rho$ on $\rbb_+$ with compact
support such that
   \begin{align} \label{cnin}
c_n= \int_{\rbb_+} x^n \D \rho(x), \quad n\in \zbb_+.
   \end{align}
   \end{lem}
   \begin{proof}
Suppose, to the contrary, that \eqref{cnin} holds for
some finite signed Borel measure $\rho$ on $\mathbb
R_+$ with compact support. Observe that
   \[
0 = (\triangle^m \boldsymbol{c})_n = \int_{\mathbb
R_+} x^n (x-1)^m \, \D\rho(x), \quad n \in \zbb_+.
   \]
Using the uniqueness of compactly supported signed
representing measures (see
\cite[Lemma~4.1]{C-J-J-S21}), we obtain
$(x-1)^m\,\D\rho(x)=0$. Hence, $\rho$ is supported on
$\{1\}$, so the sequence $\boldsymbol{c}$ is constant,
a contradiction.
    \end{proof}
   \begin{rem} \label{negj}
1) Regarding the assumption $\triangle^m
\boldsymbol{c} \equiv 0$, note the following (see
\cite[Exercise~7.2]{Dick}; see also
\cite[Proposition~2.1]{J-J-S20}):
   \begin{align*}
   \begin{minipage}{70ex}
{\em If $\boldsymbol{c}=\{c_n\}_{n=0}^\infty \subseteq
\mathbb R$ and $m\in \nbb$, then $\triangle^m
\boldsymbol{c} \equiv 0$ if and only if $c_n = p(n)$
for some real polynomial $p$ of degree at most $m-1$.}
   \end{minipage}
   \end{align*}
As a consequence, $\triangle^m \boldsymbol{c} \equiv
0$ implies $\limsup_{n\to \infty} |c_n|^{1/n} <
\infty$.

2) For any integer $m \Ge 2$ and any non-isometric
$m$-isometry $T \in \ogr{\hh}$, there exists $h \in
\hh$ such that the sequence $\{\|T^n
h\|^2\}_{n=0}^{\infty}$ satisfies the assumptions of
Lemma~\ref{trkz} (see \cite[Theorem~3.3]{J-J-S20}; see
also \cite[Theorem~4.4]{J-J-S20}).
   \end{rem}
   \section{\label{Sec.8}The  commutant criterion for subnormality}
   \begin{proof}[Proof of Theorem~\ref{cvxw}]
Denote by $G_N$ the spectral measure of $N$. Note that
the measure $G_{\phi}$ defined by
$G_{\phi}(\varDelta):=G_N (\phi^{-1}(\varDelta))$ for
all $\varDelta \in \borel{\rbb_+}$, where $\phi$ is as
in \eqref{gnms}, is the spectral measure of $|N|^2$.
Since $N$ is a minimal normal extension of $S$, we
have $\|S\|=\|N\|=\||N|\|$ (see
\cite[Corollary~II.2.17]{Co91}).

(i)$\Rightarrow$(ii) First we show that
$G_N(\cbb\setminus\dbb)=0$. By
Theorem~\ref{mnteo}~(ii-b), we have
   \allowdisplaybreaks
   \begin{align*}
0 = \mu_{|E|f}([1, \infty)) & \overset{\eqref{gnms}} =
\|G_N (\cbb \setminus \dbb) |E|f\|^2, \quad f \in
\hh_2.
   \end{align*}
This implies that for all $\varDelta \in \borel{\cbb}$
and $f \in \hh_2$,
   \begin{align} \label{gncmd}
G_N(\cbb \setminus \dbb)\left(G_N(\varDelta)|E|f
\right) =G_N(\varDelta)G_N(\cbb \setminus \dbb)|E|f=0.
   \end{align}
Since $\jd{E}=\{0\}$, $\overline{\ob{|E|}}=\hh_{2}$.
By the minimality of $N$ (see
\cite[Theorem~3]{St-SzIII}),
   \begin{align*}
\kk=\bigvee_{\varDelta \in \borel{\cbb}}
G_{N}(\varDelta)\hh_2 =\bigvee_{\varDelta \in
\borel{\cbb}} G_{N}(\varDelta)|E|\hh_2.
   \end{align*}
Together with \eqref{gncmd}, this implies that
$G_N(\cbb\setminus\dbb)=0$. Thus, by
\cite[Theorem~6.6.4]{Bi-So87},
   \begin{align} \label{gph1}
G_{|N|}([1,\infty)) = G_{\phi}([1,\infty)) =
G_N(\cbb\setminus\dbb) = 0.
   \end{align}
Therefore, $\|S\|=\||N|\|\Le 1$, and by
\cite[Theorem~6.1.3]{Bi-So87},
$1\notin\sigma_{\mathrm{p}}(|N|)$.

Now observe that
   \begin{align} \label{mhyd}
\mm := \bigvee_{\varDelta \in \borel{[0,1)}}
G_{\phi}(\varDelta) \hh_2 = \bigvee_{\varDelta \in
\borel{[0,1)}} \int_{\varDelta} \sqrt{1-x}\,
G_{\phi}(\D x)\, \hh_2.
   \end{align}
Indeed, by \cite[Theorem~6.6.3]{Bi-So87}, $\mcal$
reduces $|N|^2$ and hence, by the square root theorem,
$B:=(I-|N|^2)^{1/2}$. Since, by \eqref{gph1}, $B$ is
injective, so is $B|_{\mcal}$. Thus
   \begin{align*}
\bigvee_{\varDelta\in
\borel{[0,1)}}\int_{\varDelta}\sqrt{1-x}\, G_{\phi}(\D
x)\,\hh_2 =\bigvee_{\varDelta\in
\borel{[0,1)}}BG_{\phi}(\varDelta)\,\hh_2
=\overline{B(\mcal)}=\mcal.
   \end{align*}

A standard measure-theoretic argument applied to
Theorem~\ref{mnteo}~(ii-a) yields
   \begin{align*}
\mu_{|E| f}(\varDelta) = \int_{\varDelta}
\frac{1-x}{1-x}\,\D\mu_{|E|f}(x) \Le \int_{\varDelta}
(1-x) \,\D\mu_{f}(x), \quad \varDelta \in
\borel{[0,1)}, \, f \in \hh_2.
   \end{align*}
Together with \eqref{gnms}, this implies that
   \begin{equation} \label{eq:2}
\Big\|G_{\phi}(\varDelta)\, |E| f \Big\|^2 \Le \Big\|
\int_{\varDelta} \sqrt{1-x}\, G_{\phi}(\D x)\, f
\Big\|^2, \qquad \varDelta \in \borel{[0,1)},\, f \in
\hh_2.
   \end{equation}
We now show that for all finite sequences
$\{\varDelta_j\}_{j=1}^n \subseteq \borel{[0,1)}$ and
$\{f_j\}_{j=1}^n \subseteq \hh_2$,
   \begin{equation} \label{eq:3}
\Big\| \sum_{j=1}^n G_{\phi}(\varDelta_j)\, |E| f_j \Big\|^2 \Le
\Big\| \sum_{j=1}^n \int_{\varDelta_j} \sqrt{1-x}\, G_{\phi}(\D x)\,
f_j \Big\|^2.
   \end{equation}
Indeed, there exist pairwise disjoint sets
$\tau_1,\dots,\tau_m \in \borel{[0,1)}$ and
corresponding subsets $J_1, \ldots, J_n \subseteq
\{1,\dots,m\}$ such that
   \begin{equation} \label{eq:4}
\varDelta_j = \bigsqcup_{i \in J_j} \tau_i, \quad j=1,\dots,n.
   \end{equation}
Set $g_i=\sum_{j=1}^n \chi_{J_j}(i)\, f_j$ for $i=1,
\ldots m$. Since $G_{\phi}$ is a spectral measure, we
have
   \allowdisplaybreaks
   \begin{align*}
\Big\| \sum_{j=1}^n G_{\phi}(\varDelta_j)\, |E| f_j \Big\|^2
&\overset{\eqref{eq:4}}= \Big\| \sum_{j=1}^n \sum_{i=1}^m
\chi_{J_j}(i)\,
G_{\phi}(\tau_i)\, |E| f_j \Big\|^2
   \\
&\hspace{1ex}= \Big\| \sum_{i=1}^m G_{\phi}(\tau_i)\, |E| g_i\Big\|^2
=\sum_{i=1}^m \Big\| G_{\phi}(\tau_i)\, |E| g_i\Big\|^2
   \\
& \overset{\eqref{eq:2}} \Le \sum_{i=1}^m \Big\| G_{\phi}(\tau_i)
\int_{[0,1)} \sqrt{1-x} \, G_{\phi}(\D x)\, g_i\Big\|^2
   \\
& \hspace{1ex}= \Big\| \sum_{i=1}^m \int_{\tau_i}
\sqrt{1-x} \, G_{\phi}(\D x)\, g_i \Big\|^2
   \\
& \hspace{1ex} = \Big\| \sum_{j=1}^n \sum_{i=1}^m
\chi_{J_j}(i)\, \int_{\tau_i} \sqrt{1-x} \,
G_{\phi}(\D x)\, f_j \Big\|^2
   \\
& \overset{\eqref{eq:4}}= \Big\| \sum_{j=1}^n
\int_{\varDelta_j} \sqrt{1-x}\, G_{\phi}(\D x)\, f_j
\Big\|^2.
   \end{align*}
In view of \eqref{mhyd} and \eqref{eq:3}, there exists
a unique contraction $\tilde{R} \in \ogr{\mm}$
such~that
   \begin{align} \label{rtgph}
\tilde{R} \int_{\varDelta} \sqrt{1-x}\, G_\phi(\D x)\,
f = G_\phi(\varDelta)\, |E| f, \qquad \varDelta \in
\borel{[0,1)}, \, f \in \hh_2.
   \end{align}
Define the operator $R \in \ogr{\kk}$ by
   \begin{align*}
R = \tilde{R} \oplus 0_{\mm^{\perp}},
   \end{align*}
where $0_{\mm^{\perp}}$ is the zero operator on
$\mm^{\perp}=\kk \ominus \mm$. Then $R$ is a
contraction satisfying
   \begin{align} \label{rdpi}
R \int_{\varDelta} \sqrt{1-x}\, G_\phi(\D x)\, f =
G_\phi(\varDelta)\, |E| f, \qquad \varDelta \in
\borel{[0,1)}, \, f \in \hh_2.
   \end{align}

Now we show that $R \in \{|N|\}'$. Since $\mm$ reduces
$G_\phi$, we obtain
   \begin{align} \label{rgpt}
R\, G_\phi(\tau) = \widetilde{R}\, G_\phi(\tau)|_{\mm} \oplus
0_{\mm^\perp}, \quad \tau \in \borel{\rbb_+}.
   \end{align}
Note that, by \eqref{mhyd} and \eqref{rtgph}, we have
   \allowdisplaybreaks
   \begin{align*}
\widetilde{R}\, G_\phi(\tau) \left( \int_{\varDelta} \sqrt{1-x}\,
G_\phi(\D x)\, f \right) &= \widetilde{R}
\int_{\tau \cap \varDelta} \sqrt{1-x}\, G_\phi(\D x)\, f
   \\
& = G_\phi(\tau) \bigl(G_\phi(\varDelta)|E|f\bigr)
   \\
& = G_\phi(\tau)\, \widetilde{R} \left(
\int_{\varDelta} \sqrt{1-x}\, G_\phi(\D x)\, f
\right),
   \end{align*}
for all $\tau \in \borel{\rbb_+}$, $\varDelta \in
\borel{[0,1)}$, and $f \in \hh_2$. This, together with
\eqref{mhyd}, yields
   \[
\widetilde{R}\, G_\phi(\tau)|_{\mm} = G_\phi(\tau)|_{\mm}\,
\widetilde{R}, \quad \tau \in \borel{\rbb_+}.
   \]
Hence, by \eqref{rgpt} and
\cite[Theorem~6.6.3]{Bi-So87}, $R \in
\bigl\{G_\phi(\tau)\colon \tau \in
\borel{\rbb_+}\bigr\}'= \{|N|\}^{\prime}$.

Using the functional calculus and combining
\eqref{gph1} with \eqref{rdpi}, we obtain
   \begin{align*}
R (I - |N|^2)^{1/2}|_{\hh_2} f= R \int_{[0,1)} \sqrt{1-x}\, G_\phi(\D
x)\, f = |E| f, \quad f \in \hh_2,
   \end{align*}
which yields \eqref{weqsa}. By \eqref{gph1}, $0 \Le
|N|^2 \Le I$, so $(I - |N|^2)^{1/2}$ is a contraction.
Since $\|R\| \Le 1$, \eqref{weqsa} gives $\|E\| \Le
1$. This proves the ``moreover'' part.

(ii)$\Rightarrow$(i) Since $\bigl\||N|^2\bigr\| =
\|S\|^2 \Le 1$, we obtain $G_{\phi}((1, \infty))=0$.
As $\|R\| \Le 1$ and $R \in \{|N|^2\}' =
\bigl\{G_{\phi}(\varDelta)\colon \varDelta \in
\borel{\rbb_+}\bigr\}'$, we deduce from \eqref{gnms}
and \eqref{weqsa} that
   \begin{align*}
\mu_{|E|f}(\varDelta) & = \| G_\phi(\varDelta) |E| f
\|^2
   \\
& = \|G_\phi(\varDelta) R \bigl(I - |N|^2\bigr)^{1/2}
f \|^2
   \\
& = \Bigl \|R \int_{\varDelta} \sqrt{1-x} \, G_\phi(\D
x) f\Bigl\|^2
   \\
& \Le \int_{\varDelta} (1-x)\, \D \mu_f(x), \qquad
\varDelta \in \borel{[0,1]}, \, f \in \hh_2.
   \end{align*}
Integrating the function $(1-x)^{-1}$ relative to
$\mu_{|E|f}$ and $(1-x)\, \D \mu_f(x)$, we obtain the
inequality (ii-a) of Theorem~\ref{mnteo}. Taking
$\varDelta=\{1\}$ yields the equality (ii-b) of this
theorem, which implies the subnormality of $T$.
   \end{proof}
   \begin{exa}
The subspace $\lcal$ need not reduce an arbitrary
contraction $R\in\{|N|\}'$ satisfying
$R(I-|N|^2)^{1/2}f=|E|f$ for all $f\in\hh_2$. Indeed,
fix $a\in(0,1)$. Let $\{e_n\}_{n\in\mathbb Z}$ be the
standard orthonormal basis of $\kk=\ell^2(\mathbb Z)$,
and let $U\in\ogr{\kk}$ be the unitary bilateral
shift, i.e., $Ue_n=e_{n+1}$ for all $n\in\zbb$. Set
$\hh_2=\bigvee\nolimits_{n\Ge0}\cbb\cdot e_n$ and
$N=aU$. Then $N$ is a minimal normal extension of
$S:=N|_{\hh_2}$, and \mbox{$|N|=aI_{\kk}$}.~Thus
\[
\lcal:=\bigvee\nolimits_{n\Ge0}|N|^n\hh_2=\hh_2, \quad
B:=(I-|N|^2)^{1/2}=\sqrt{1-a^2}\,I_{\kk}.
\]
Take $|E|=a\sqrt{1-a^2}\,I_{\hh_2}$, with positive
lift $\widehat{|E|}=a\sqrt{1-a^2}\,I_{\kk}$. Define
\[
F=e_0\otimes e_{-1},\quad R=aI_{\kk}+(1-a)F.
\]
Here $(u\otimes v)f=\langle f,v\rangle u$ for all
$u,v,f\in\kk$. Clearly, $\|R\|\Le a+(1-a)=1$ and
$R\in\{|N|\}'$. Moreover, $F|_{\hh_2}=0$, and hence
\[
RBf=a\sqrt{1-a^2}\,f=|E|f, \quad f\in\hh_2.
\]
Thus $\lcal$ is invariant under $R$, but it does not
reduce $R$, because
\[ \tag*{$\diamondsuit$}
R^*e_0=ae_0+(1-a)e_{-1}\notin\lcal.
\]
   \end{exa}
   \begin{proof}[Proof of Theorem~\ref{wkwsc}]
Proceeding as in the first paragraph of the proof of
Lemma~\ref{snqm} and using Appendix~\ref{Subs.11.3},
we can assume without loss of generality that $M_1=N$,
$M_2=\widehat{|E|}$, $Q_1=N|_{\lcal}$, and
$Q_2=\widehat{|E|}|_{\lcal}$, where $N\in \ogr{\kk}$
is a minimal normal extension of $S$, $\widehat{|E|}$
is the lift of $|E|$ relative to $\{N\}^{\prime}$, and
   \begin{align*}
\lcal := \bigvee_{n\Ge 0} |N|^n \hh_2.
   \end{align*}
Since $|E|$ is injective, so is $\widehat{|E|}$ (see
\cite[Lemma~3.2]{C-J-J-S24}). By the minimality of
$N$,
   \begin{align*}
\||N|\|=\|N\|=\|S\|.
   \end{align*}

First we show that
   \begin{align} \label{ixie}
   \begin{minipage}{65ex}
{\em $T$ is subnormal if and only if $(|S|,|E|)$ is a
spherical contraction.}
   \end{minipage}
   \end{align}

Suppose that $T$ is subnormal. By
Theorem~\ref{cvxw}~(ii), $\|S\|\Le 1$ and there exists
a contraction $R \in \left\{|N|\right\}'$ satisfying
\eqref{weqsa}. Setting $B:=(I-|N|^2)^{1/2}$ yields
   \begin{align*}
\is{|E|^2 f}{f} & = \|R B f\|^2
   \\
& \Le \|B f\|^2
   \\
& = \|f\|^2 - \|N f\|^2
   \\
& = \|f\|^2 - \is{|S|^2 f}{f}, \quad f \in \hh_2,
   \end{align*}
which implies that $|E|^2 + |S|^2 \Le I$.

Suppose now that $(|S|,|E|)$ is a spherical
contraction. By Lemmas~\ref{hcps} and~\ref{snqm} (the
second applied to $A=|E|$), the pair
$(|N|,\widehat{|E|})$ is a spherical contraction. In
particular, $\|S\| = \||N|\| \Le1$, so
$B:=(I-|N|^2)^{1/2}$ is well defined. Moreover,
   \begin{align*}
\|\widehat{|E|}f\|^2 \Le\|f\|^2-\is{|N|^2f}{f}
=\|Bf\|^2, \quad f\in\kk.
   \end{align*}
Hence there exists a contraction $R\in\ogr{\kk}$ such
that $R B = \widehat{|E|}$. Since $\widehat{|E|}$ is
injective, the equality implies that $B$ is injective.
Moreover, $\widehat{|E|}\in\{N\}'$ commutes with $B$
by the square root theorem. Applying Lemma~\ref{comx}
to $A=\widehat{|E|}$, we obtain $R\in\{B\}'$. Since
$|N|=(I-B^2)^{1/2}$, another application of the square
root theorem yields $R\in\{|N|\}'$. Restricting $R B =
\widehat{|E|}$ to $\hh_2$ gives \eqref{weqsa}, because
$\widehat{|E|}|_{\hh_2}=|E|$. Theorem~\ref{cvxw} now
implies that $T$ is subnormal.

It follows from \eqref{ixie} and the identity $T^*T =
I_{\hh_1}\oplus\bigl(|S|^2+|E|^2\bigr)$ that $T$ is
subnormal if and only if $T$ is a contraction.

The ``moreover'' part follows directly from
\eqref{ixie} and Lemmas~\ref{hcps} and~\ref{snqm}
(also use the equality $|N|_{\lcal}| = |N||_{\lcal}$).
   \end{proof}
   \section{\label{Sec.9}Subnormality: the general case}
We begin by reducing the subnormality problem for
$B$-operators of class $\mathscr{S}$ to the case where
$E$ is injective. This reduction is inspired by
\cite[Theorem~1.8]{C-J-J-S24}.
   \begin{lem} \label{th2z}
Let $T=\left[\begin{smallmatrix} V & E \\ 0 & S
\end{smallmatrix}\right]\in\ogr{\hh_1\oplus\hh_2}$
be a $B$-operator with $E\neq 0$. Then the following
statements hold\textup{:}
   \begin{enumerate}
   \item[\textup{(i)}]
the subspaces $\hh_{2,0}:=\overline{\ob{|E|}}$ and
$\hh_{2,1}:=\jd{E}$ reduce $S$, with
$\hh_{2,0}\neq\{0\}$, and determine the orthogonal
decompositions
   \begin{align*}
\hh_2=\hh_{2,0}\oplus\hh_{2,1},\qquad S=S_0\oplus S_1,
   \end{align*}
where $S_j:=\restr{S}{\hh_{2,j}}$ for $j=0,1$,
   \item[\textup{(ii)}]
$|E_0|= \restr{|E|}{\hh_{2,0}}$,
$\overline{\ob{|E_0|}} = \hh_{2,0}$, and $E_1=0$,
where
   \begin{align*}
E_0:=\restr{E}{\hh_{2,0}}\in \ogr{\hh_{2,0},\hh_1}
\text{ and } E_1:=\restr{E}{\hh_{2,1}} \in
\ogr{\hh_{2,1},\hh_1},
   \end{align*}
   \item[\textup{(iii)}]
$T_0:=\smat{V&E_0\\0&S_0}\in\ogr{\hh_1\oplus
\hh_{2,0}}$ is a $B$-operator with $E_0$ injective,
   \item[\textup{(iv)}]
if $T$ is of class $\mathscr{S}$, then so is $T_0$,
   \item[\textup{(v)}]
$T=T_0\oplus S_1$ with respect to the orthogonal
decomposition
   \begin{align} \label{myptr}
\hh_1\oplus\hh_2=(\hh_1\oplus\hh_{2,0})\oplus
\hh_{2,1},
   \end{align}
   \item[\textup{(vi)}]
$\sigma(T)=\sigma(T_0)\cup\sigma(S_1)$,
   \item[\textup{(vii)}]
$T$ is subnormal if and only if $T_0$ and $S_1$ are
subnormal.
   \end{enumerate}
Moreover, if $T$ is of class $\mathscr{S}$, and is
subnormal, then $S_0$ is a contraction.
   \end{lem}
   \begin{proof}
Clearly, the subspace $\hh_{2,0}$ is nonzero. Since
$S$ commutes with the positive operator $|E|$, the
subspaces $\hh_{2,0}=\overline{\ob{|E|}}$ and
$\hh_{2,1}=\jd{E}$ reduce both $S$ and $|E|$. Hence
$|E_0|=\restr{|E|}{\hh_{2,0}}$, and therefore $S_0$
commutes with $|E_0|$. Moreover,
\[
\jd{E_0}=\{0\}, \quad \overline{\ob{|E_0|}}=\hh_{2,0}
\;\; \text{and} \;\; \ob{E_0}=\ob{E}\perp\ob{V}.
\]
These properties show that $T_0$ is a $B$-operator
with $E_0$ injective.

With respect to the orthogonal decomposition
$\hh_2=\hh_{2,0}\oplus\hh_{2,1}$, we have
\[
S=S_0\oplus S_1 \;\; \text{and} \;\; |E|=|E_0|\oplus
0.
\]
This proves \textup{(i)}-\textup{(iii)}. Statement
\textup{(iv)} is now obvious.

Observe that the orthogonal decomposition
\eqref{myptr} holds, with all direct sums understood
as internal orthogonal direct sums. A routine
verification shows that, with respect to this
decomposition,
\[
T=T_0\oplus S_1,
\]
which is precisely \textup{(v)}. The remaining
statements \textup{(vi)}-\textup{(vii)} follow from
\textup{(v)}.

The ``moreover'' part is a consequence of the above
arguments and Theorem~\ref{cvxw} applied to the
$B$-operator $T_0$.
   \end{proof}
   \begin{proof}[Proof of Theorem~\ref{aswq}]
Proceeding as in the first paragraph of the proof of
Lemma~\ref{snqm} and using Appendix~\ref{Subs.11.3},
we can assume without loss of generality that $M_1=N$,
$M_2=\widehat{|E|}$, $Q_1=N|_{\lcal}$, and
$Q_2=\widehat{|E|}|_{\lcal}$, where $N\in\ogr{\kk}$ is
a minimal normal extension of $S$, $\widehat{|E|}$ is
the lift of $|E|$ relative to $\{N\}^{\prime}$, and
   \begin{align*}
\lcal:=\bigvee_{n\Ge0}|N|^n\hh_2.
   \end{align*}

By Lemma~\ref{th2z}, we have the following orthogonal
decompositions:
   \begin{align*}
\hh_2=\hh_{2,0}\oplus\hh_{2,1},\quad S=S_0\oplus
S_1,\quad |E|=|E_0|\oplus0.
   \end{align*}
Thus
   \begin{align} \label{dwqa}
(S,|E|)=(S_0,|E_0|)\oplus(S_1,0).
   \end{align}
By \cite[Lemma~3.4(i)]{C-J-J-S24}, for $j=0,1$, the
subspace $\kk_j:=\bigvee_{m\Ge0}N^{*m}\hh_{2,j}$
reduces $N$, and $N_j:=\restr{N}{\kk_j}$ is a minimal
normal extension of $S_j$. Moreover,
   \begin{align} \label{cor25-normal-decomposition}
\kk=\kk_0\oplus\kk_1 \quad\text{and}\quad N=N_0\oplus
N_1.
   \end{align}
Clearly, $|N_j|=|N||_{\kk_j}$ for $j=0,1$. Recall that
$E_1$ is the zero operator from $\hh_{2,1}$ into
$\hh_1$. By \cite[Lemma~3.4(ii)]{C-J-J-S24}, for
$j=0,1$, the subspace $\kk_j$ reduces $\widehat{|E|}$,
and $\widehat{|E|}|_{\kk_j}$ is the lift of $|E_j|$
relative to $\{N_j\}^{\prime}$, that is,
$\widehat{|E|}|_{\kk_j}=\widehat{|E_j|}$. Hence
$\widehat{|E|}=\widehat{|E_0|}\oplus0$, and
   \begin{align} \label{cor25-lift-decomposition}
(N,\widehat{|E|}) =(N_0,\widehat{|E_0|})\oplus(N_1,0).
   \end{align}
Lemma~\ref{liom} shows that $(N,\widehat{|E|})$ and
$(N_0,\widehat{|E_0|})$ are minimal normal extensions
of $(S,|E|)$ and $(S_0,|E_0|)$, respectively. For
$j=0,1$, the subspace $\lcal_j$ of $\kk_j$ defined by
   \begin{align} \label{mycp}
\lcal_j :=\bigvee_{m,n\Ge0}|N_j|^m
\widehat{|E_j|}^{\,n}\hh_{2,j}
=\bigvee_{m\Ge0}|N|^m\hh_{2,j}
   \end{align}
is invariant under $N_j$ and $\widehat{|E_j|}$, and
the pair
   \begin{align} \label{vypx}
(N_j|_{\lcal_j},\widehat{|E_j|}|_{\lcal_j})
\overset{\eqref{cor25-lift-decomposition}}{=}
(N|_{\lcal_j},\widehat{|E|}|_{\lcal_j})
   \end{align}
is a minimal quasinormal extension of $(S_j,|E_j|)$.
By \eqref{mycp} and
\eqref{cor25-normal-decomposition},
\mbox{$\lcal=\lcal_0\oplus\lcal_1$},~so
   \begin{align} \label{przcs}
(N|_{\lcal},\widehat{|E|}|_{\lcal})
=(N_0|_{\lcal_0},\widehat{|E_0|}|_{\lcal_0})
\oplus(N_1|_{\lcal_1},0).
   \end{align}
Since $E_0$ is injective, so is $|E_0|$. By
\cite[Lemma~3.2]{C-J-J-S24}, the lift
$\widehat{|E_0|}$ is therefore injective, as is its
restriction $\widehat{|E_0|}|_{\lcal_0}$. Thus the
pairs $(S_0,|E_0|)$, $(N_0,\widehat{|E_0|})$, and
$(N_0|_{\lcal_0},\widehat{|E_0|}|_{\lcal_0})$ satisfy
the assumptions of Lemma~\ref{cxzq}.

Suppose that $T$ is subnormal. By
Proposition~\ref{13posr}, $T_0$ is subnormal. Hence
Theorem~\ref{wkwsc} yields
$\sigma(S_0,|E_0|)\subseteq\sbd$. Using \eqref{ossp},
\eqref{dwqa}, and the projection property of the
Taylor spectrum (see \cite[p.\ 21]{Cur88}), we obtain
   \begin{align*}
\sigma(S,|E|) =\sigma(S_0,|E_0|)\cup\sigma(S_1,0)
\subseteq\varXi.
   \end{align*}

Conversely, suppose that
$\sigma(S,|E|)\subseteq\varXi$. It follows from
\eqref{ossp} and \eqref{dwqa} that
$\sigma(S_0,|E_0|)\subseteq\varXi$. Lemma~\ref{cxzq}
then gives $\sigma(S_0,|E_0|)\subseteq\sbd$. By
Theorem~\ref{wkwsc} and Proposition~\ref{13posr}, $T$
is subnormal.

Repeating the above reasoning for the pairs
$(N_0,\widehat{|E_0|})$ and
$(N_0|_{\lcal_0},\widehat{|E_0|}|_{\lcal_0})$ and
using the decompositions
\eqref{cor25-lift-decomposition} and \eqref{przcs}, we
obtain the characterizations of subnormality of $T$
given by (iii) and (v), respectively. Finally, the
decompositions \eqref{dwqa},
\eqref{cor25-lift-decomposition}, and \eqref{przcs}
remain valid when the first component of each pair is
replaced by its modulus. The same argument, now using
the positive case of Lemma~\ref{cxzq}, yields the
characterizations in (ii), (iv), and (vi),
respectively.
   \end{proof}
In the following corollary, $(0,1)$ is understood as a
point of $\cbb^2$.
   \begin{cor}
   Let $T = \left[\begin{smallmatrix} V & E \\ 0 & S
\end{smallmatrix}\right] \in \ogr{\hh_1 \oplus \hh_2}$
be a $B$-operator of class $\mathscr{S}$ with $E \neq
0$, and let $\nbold=(N_1,N_2)$ be a minimal normal
extension of $(S,|E|)$. Assume that $|E|$ is an
orthogonal projection. Then the following conditions
are equivalent{\em :}
   \begin{enumerate}
   \item[(i)] $T$ is subnormal,
   \item[(ii)] $S|E|=0$ $($equivalently,
   $|S||E|=0$$)$,
   \item[(iii)] $\sigma(|S|, |E|) \subseteq \{(0,
1)\} \cup \big(\mathbb R_+ \times \{0\}\big)$,
   \item[(iv)] $\sigma(S, |E|) \subseteq \{(0,
1)\} \cup \big(\mathbb C \times \{0\}\big)$,
   \item[(v)] $\sigma(\nbold) \subseteq
\varXi$.
   \end{enumerate}
   \end{cor}
   \begin{proof}
(i)$\Rightarrow$(ii) Let $\lcal_0$ be as in
\eqref{mycp}. First, note that
   \begin{align*}
\widehat{|E|}(|N|^j |E|h) = |N|^j |E|^2h = |N|^j |E|h,
\qquad h \in \hh_2, \, j \Ge 0.
   \end{align*}
It follows that $\widehat{|E|}|_{\lcal_0}=I$. Since
$(N|_{\lcal_0},\widehat{|E|}|_{\lcal_0})$ is a minimal
quasinormal extension of $(S_0,|E_0|)$ (see
\eqref{vypx}), Theorem~\ref{13posr} implies that the
pair $\bigl(\bigl|N|_{\lcal_0}\bigr|,
\widehat{|E|}|_{\lcal_0}\bigr)$ is a spherical
contraction. Thus $\bigl|N|_{\lcal_0}\bigr|^2+I \Le
I$, which forces $N|_{\lcal_0}=0$. Since $\ob{|E|}
\subseteq \lcal_0$, we conclude that
   \begin{align*}
S|E|h=N|E|h=0, \quad h\in \hh_2,
   \end{align*}
which implies (ii).

(ii)$\Rightarrow$(iv) Set $P:=|E|$ and note that
$\sigma(P) \subseteq \{0,1\}$. By the spectral mapping
theorem and projection property for the Taylor
spectrum (see \cite[Theorem~4.8]{Tay70-c}; see also
\cite[p.\ 21]{Cur88}), the equality $SP=0$ yields
   \begin{align*}
\sigma(S,P) &\subseteq \big\{(z,w)\in \mathbb C^2
\colon zw=0, \, z\in \sigma(S), \, w\in
\sigma(P)\big\}
   \\
&\subseteq \big\{(z,w)\in \mathbb C^2 \colon zw=0, \,
w\in \{0,1\}\big\}
   \\
&= \{(0,1)\}\cup \big(\mathbb C\times \{0\}\big).
   \end{align*}

(iv)$\Rightarrow$(v) This follows from the spectral
inclusion theorem for subnormal pairs (see
\eqref{pytr}).

(v)$\Rightarrow$(i) This is a direct consequence of
Theorem~\ref{aswq}.

(ii)$\Rightarrow$(iii) The proof is analogous to that
of (ii)$\Rightarrow$(iv).

(iii)$\Rightarrow$(i) This follows again from
Theorem~\ref{aswq}.
   \end{proof}
   \appendix
   \section{\label{Sec.A}Joint quasinormality}
In this section, we discuss the basic properties of quasinormal pairs
of commuting operators, along the lines of the theory for single
operators developed in \cite{Em81}. For this purpose, we found it
necessary to choose among several definitions of quasinormality
appearing in \cite{Em86,Gl06} (see also \cite{Ath2015}); our choice
is \cite[Definition~1.1(b)]{Gl06}. The results established below are
not available elsewhere in the literature and extend verbatim to
$n$-tuples.
    \subsection{Quasinormal pairs}
We say that a pair $\tbold=(T_1,T_2)\in \ogr{\hh}^2$ of commuting
operators is {\em quasinormal} if $T_iT_j^*T_j=T_j^*T_jT_i$ for
$i,j=1,2$, or equivalently, if $T_i|T_j|=|T_j|T_i$ for $i,j=1,2$.
Clearly, every quasinormal pair consists of commuting quasinormal
operators. By \cite[Proposition~2.1]{Ath2015}, every quasinormal pair
is subnormal. This leads to the following result (compare with
\cite[Theorem~2.9]{Cur-Lee-Yoon20}).
   \begin{pro}  \label{sdaq}
Let $\tbold=(T_1,T_2)\in \ogr{\hh}^2$ be a pair of commuting
operators. Then the following conditions are equivalent{\em :}
   \begin{enumerate}
   \item[(i)] $\tbold$ is quasinormal,
   \item[(ii)] $\tbold$ is subnormal and both $T_1$ and $T_2$ are
quasinormal.
   \end{enumerate}
   \end{pro}
   \begin{proof}
By the above discussion, (i) implies (ii). If (ii)
holds and $\nbold=(N_1,N_2)\in \ogr{\kk}^2$ is a
normal extension of $\tbold$, then, by \eqref{campe},
$|N_j|(\hh) \subseteq \hh$ and thus
$|T_j|=|N_j||_{\hh}$ ($j=1,2$), which implies that
$\tbold$ is quasinormal.
   \end{proof}
   \subsection{\label{Subs.11.2}Minimal quasinormal extensions}
We say that a quasinormal extension $\qbold=(Q_1,Q_2)\in
\ogr{\lcal}^2$ of a pair $\tbold=(T_1,T_2)\in \ogr{\hh}^2$ of
commuting operators is {\em minimal} if $\lcal$ is the only closed
subspace $\tilde{\lcal}\subseteq\lcal$ containing $\hh$, invariant
under $Q_1$ and $Q_2$, and such that the pair
$(Q_1|_{\tilde{\lcal}},Q_2|_{\tilde{\lcal}})$ is quasinormal.

   Let $\sbold=(S_1,S_2)\in \ogr{\hh}^2$ be a
subnormal pair, and let $\nbold=(N_1,N_2) \in
\ogr{\kk}^2$ be its normal extension. Since, by the
Fuglede theorem, the operators $N_1$, $N_1^*$, $N_2$
and $N_2^*$ commute, the space
   \begin{align*}
\lcal:=\bigvee_{m,n\Ge 0} |N_1|^{m} |N_2|^{n} \hh
   \end{align*}
is invariant under $N_j$ and $|N_j|$ for $j=1,2$.
Hence, by \eqref{campe}, the operator $Q_j =
N_j|_{\lcal}$ is quasinormal for $j=1,2$. Clearly,
$Q_1$ and $Q_2$ commute. Since $|Q_j|=|N_j||_{\lcal}$
for $j=1,2$, the pair $\qbold:=(Q_1,Q_2)$ is
quasinormal (this can also be deduced from
Proposition~\ref{sdaq}). We show that the pair
$\qbold$ is a minimal quasinormal extension of
$\sbold$. Indeed, suppose that a closed subspace
$\tilde{\lcal}$ satisfies $\hh \subseteq \tilde{\lcal}
\subseteq \lcal$ and $Q_j(\tilde{\lcal}) \subseteq
\tilde{\lcal}$ for $j=1,2$, and that the pair
$(\tilde{Q_1}, \tilde{Q_2})$ is quasinormal, where
$\tilde{Q}_j := Q_j|_{\tilde{\lcal}} =
N_j|_{\tilde{\lcal}}$ for $j=1,2$. Applying
\eqref{campe} to the quasinormal operators
$\tilde{Q}_1$ and $\tilde{Q}_2$, we obtain
$N_j(\tilde{\lcal}) \subseteq \tilde{\lcal}$ and
$|N_j|(\tilde{\lcal}) \subseteq \tilde{\lcal}$ for
$j=1,2$. This implies that
   \begin{align*}
\lcal = \bigvee_{m,n\Ge 0} |N_1|^{m} |N_2|^{n} \hh
\subseteq \bigvee_{m,n \Ge 0} |N_1|^m |N_2|^n
\tilde{\lcal} \subseteq \tilde{\lcal} \subseteq \lcal,
   \end{align*}
which shows that $\lcal = \tilde{\lcal}$. If $\nbold$
is a minimal normal extension of $\sbold$, then
   \begin{align*}
\kk \overset{\eqref{mnex}}= \bigvee_{m,n \Ge 0}
N_1^{*m} N_2^{*n} \hh \subseteq \bigvee_{m,n \Ge 0}
N_1^{*m} N_2^{*n} \lcal \subseteq \kk,
   \end{align*}
and hence $\nbold$ is a minimal normal extension of
$\qbold$.
   \subsection{\label{Subs.11.3}$\hh$-unitary equivalence}
We now show that minimal quasinormal extensions of a
subnormal pair $\sbold=(S_1, S_2)\in \ogr{\hh}^2$ are
$\hh$-unitarily equivalent (with the same meaning as
for minimal normal extensions). Let $\qbold = (Q_1,
Q_2) \in \ogr{\lcal}^2$ and $\qbold^{\prime} =
(Q_1^{\prime}, Q_2^{\prime}) \in
\ogr{\lcal^{\prime}}^2$ be arbitrary minimal
quasinormal extensions of $\sbold$. Since quasinormal
pairs are subnormal, there exist normal extensions
$\nbold=(N_1,N_2)\in \ogr{\kk}^2$ and
$\nbold^{\prime}=(N_1^{\prime},N_2^{\prime})\in
\ogr{\kk^{\prime}}^2$ of $\qbold$ and
$\qbold^{\prime}$, respectively. Clearly, these are
also normal extensions of $\sbold$. Using
\eqref{campe}, we deduce that $\lcal$ is invariant
under $N_j$ and $|N_j|$, and consequently
$|Q_j|=|N_j||_{\lcal}$ for $j=1,2$. In particular,
$|N_j|(\lcal) \subseteq \lcal$ for $j=1,2$. This
implies that
   \begin{align*}
\lcal_0:=\bigvee_{m,n \Ge 0} |N_1|^{2m} |N_2|^{2n} \hh
\subseteq \bigvee_{m,n \Ge 0} |N_1|^m |N_2|^n \lcal
\subseteq \lcal.
   \end{align*}
Clearly, $\lcal_0$ is invariant under $|N_1|^2$ and
$|N_2|^2$, and hence, by the square root theorem, it
is invariant under $|N_1|$ and $|N_2|$. Since
$\lcal_0$ is also invariant under $N_1$ and $N_2$, we
infer from \eqref{campe} that $|N_j|_{\lcal_0}|=
|N_j||_{\lcal_0}$, and thus
$(Q_1|_{\lcal_0},Q_2|_{\lcal_0})$ is a quasinormal
pair. The minimality of $\qbold$ implies that
$\lcal_0=\lcal$. Therefore, we obtain
   \begin{align} \label{Minx}
\lcal=\bigvee_{m,n \Ge 0} |N_1|^{2m} |N_2|^{2n} \hh.
   \end{align}
A similar argument shows that
$\lcal^{\prime}=\bigvee_{m,n \Ge 0}
|N_1^{\prime}|^{2m} |N_2^{\prime}|^{2n} \hh$.
Therefore, for all $f,g \in \hh$ and all $m,n,k,l \in
\zbb_+$,
   \allowdisplaybreaks
   \begin{align*}
\langle |N_1|^{2m} |N_2|^{2n} f, |N_1|^{2k} |N_2|^{2l}
g\rangle & = \langle S_1^{m+k} S_2^{n+l} f, S_1^{m+k}
S_2^{n+l} g\rangle
   \\
& = \langle |N_1^{\prime}|^{2m} |N_2^{\prime}|^{2n} f,
|N_1^{\prime}|^{2k} |N_2^{\prime}|^{2l} g\rangle.
   \end{align*}
This implies that there exists a unique $\hh$-unitary
isomorphism $U\colon \lcal \to \lcal^{\prime}$
such~that
   \begin{align*}
U\bigl(|N_1|^{2m} |N_2|^{2n} f\bigr) =
|N_1^{\prime}|^{2m} |N_2^{\prime}|^{2n} f, \qquad f\in
\hh, \, m,n \in \zbb_+.
   \end{align*}
Since, for all $f \in \hh$ and all $m,n \in \zbb_+$,
   \begin{align*}
UN_j \bigl(|N_1|^{2m} |N_2|^{2n} f\bigr) =
N_j^{\prime} U \bigl(|N_1|^{2m} |N_2|^{2n} f\bigr),
\quad j=1,2,
   \end{align*}
we conclude that $\qbold$ and $\qbold^{\prime}$ are
$\hh$-unitarily equivalent.
   \subsection{Nesting minimality}
   Our final observation is that if $\nbold$ is a minimal
normal extension of a minimal quasinormal extension
$\qbold$ of $\sbold$, then $\nbold$ is also a minimal
normal extension of $\sbold$. Indeed, for all $h \in
\hh$ and all $m,n,k,l \in \zbb_+$,
   \begin{align*}
N_1^{*m} N_2^{*n} |N_1|^{2k} |N_2|^{2l} h =
N_1^{*(m+k)} N_2^{*(n+l)} N_1^k N_2^l h = N_1^{*(m+k)}
N_2^{*(n+l)} S_1^k S_2^l h,
   \end{align*}
and hence
   \begin{align*}
\kk = \bigvee_{m,n \Ge 0} N_1^{*m} N_2^{*n} \lcal
\overset{\eqref{Minx}}{\subseteq} \bigvee_{m,n \Ge 0}
N_1^{*m} N_2^{*n} \hh \subseteq \kk,
   \end{align*}
which proves that $\nbold$ is a minimal normal
extension of $\sbold$.
   
   \end{document}